\documentclass{article}
\usepackage{amsmath,amsthm,bm,colortbl,graphicx,mathrsfs,afterpage,amssymb,url}
\pdfoutput=1
\theoremstyle{plain}
\newtheorem{theorem}{Theorem}
\newtheorem{lemma}[theorem]{Lemma}
\newtheorem{proposition}[theorem]{Proposition}
\newtheorem{conjecture}[theorem]{Conjecture}
\theoremstyle{remark} 
\newtheorem*{remark}{Remark} 
\def\T{\textsf{T}} 
\def\eps{\varepsilon} 
\allowdisplaybreaks
\makeatother
\title{A discrete dynamical system for the greedy strategy at collective Parrondo games}
\author{S. N. Ethier\thanks{Department of Mathematics, University of Utah, 155 S. 1400 E., Salt Lake City, UT 84112, USA. e-mail: ethier@math.utah.edu.}\; and Jiyeon Lee\thanks{Department of Statistics, Yeungnam University, 214-1 Daedong, Kyeongsan, Kyeongbuk 712-749, South Korea.  e-mail: leejy@yu.ac.kr.}}
\date{}

\begin{document}
\maketitle

\begin{abstract}
We consider a collective version of Parrondo's games with probabilities parametrized by $\rho\in(0,1)$ in which a fraction $\phi\in(0,1]$ of an infinite number of players collectively choose and individually play at each turn the game that yields the maximum average profit at that turn.  Din\'is and Parrondo [{\em Optimal strategies in collective Parrondo games}, Europhys. Lett. 63 (2003), pp. 319--325] and Van den Broeck and Cleuren [{\em Parrondo games with strategy}, Proceedings of the SPIE 5471 (2004), pp. 109--118] studied the asymptotic behavior of this greedy strategy, which corresponds to a piecewise-linear discrete dynamical system in a subset of the plane, for $\rho=1/3$ and three choices of $\phi$.  We study  its asymptotic behavior for all $(\rho,\phi)\in(0,1)\times(0,1]$, finding that there is a globally asymptotically stable equilibrium if $\phi\le2/3$ and, typically, a unique (asymptotically stable) limit cycle if $\phi>2/3$ (``typically'' because there are rare cases with two limit cycles).  Asymptotic stability results for $\phi>2/3$ are partly conjectural.
\medskip\par

\end{abstract}

\newpage

\section{Introduction}
\label{intro}
The Parrondo effect, in which there is a reversal in direction in some system parameter when two similar dynamics are combined, is the result of an underlying nonlinearity.  It was first described by J. M. R. Parrondo in 1996 in the context of games of chance:  He showed that it is possible to combine two losing games to produce a winning one.  His idea has inspired research in such diverse areas as chemistry \cite{O}, evolutionary biology \cite{XPYX}, population genetics \cite{R}, finance \cite{S}, reliability theory \cite{D}, chaos \cite{APR}, fractals \cite{APA}, epistemology \cite{St}, quantum mechanics \cite{FA}, and probability theory \cite{Py}.  In the present paper, we analyze a discrete dynamical system, introduced by Din\'is and Parrondo \cite{DP}, that models the short-range optimization, or greedy, strategy at collective Parrondo games.  Our analysis gives conditions under which the Parrondo effect is present.

Let us describe Parrondo's original example.  The so-called capital-dependent Parrondo games consist of two games, $A$ and $B$. In game $A$, the player wins one unit with probability $1/2-\eps$, where $\eps\ge0$ is a small bias parameter, and loses one unit otherwise. Game $B$ is played with two coins, the one tossed depending on the current capital of the player: If the player's current capital is divisible by 3, he tosses a ``bad'' coin with probability of heads $p_0 := 1/10 -\eps$, otherwise he tosses a ``good'' coin with probability of heads $p_1:=3/4-\eps$; he wins one unit with heads and loses one unit with tails.  It can be shown that when $\eps=0$, both games $A$ and $B$ are fair, hence losing when $\eps>0$.  However, the random mixture $\gamma A+(1-\gamma)B$, in which game $A$ is played with probability $\gamma\in(0,1)$ and game $B$ is played with probability $1-\gamma$, is winning for $\eps\ge0$ sufficiently small. Furthermore, the non-random pattern $[r,s]$, denoting $r$ plays of game $A$ followed by $s$ plays of game $B$ (repeated ad infinitum), is also winning for $\eps\ge0$ sufficiently small, except when $r=s=1$. In summary, two losing (or fair) games can be combined, by random mixture or nonrandom alternation, to create a winning game.  This was the original form of \textit{Parrondo's paradox}.  See \cite{HA, PD, E, A} for survey articles.

Din\'is and Parrondo \cite{DP} formulated a modification of the capital-dependent Parrondo games in which a fraction of an infinite number of players collectively choose and individually play the same game at each turn. They found that, in certain cases, by choosing the game that yields the maximum average profit at each turn, the surprising result is systematic losses (if $\eps>0$), whereas a random or nonrandom sequence of choices yields a steady increase in average profit.  Van den Broeck and Cleuren \cite{VC} considered also the case of a finite number of players.  They evaluated the expected profit for this greedy strategy as a function of the number of players and proved that the strategy is optimal when the number of players is one or two but suboptimal when it is three or infinite.

In this paper we consider only the case of infinitely many players, and we adopt the one-parameter family of capital-dependent Parrondo games of Ethier and Lee \cite{EL} given by $p_0:=\rho^2/(1+\rho^2)-\eps$ and $p_1:=1/(1+\rho)-\eps$ with $\rho>0$; the original Parrondo game $B$, assumed by Din\'is and Parrondo \cite{DP} and Van den Broeck and Cleuren \cite{VC}, corresponds to $\rho=1/3$.  In order to focus on the case in which two fair games produce a winning
game, we assume that $\rho\in(0,1)$ for game $B$ and $\eps=0$ for both games.  The fraction of players who play at each turn, denoted here by $\phi$ (as in \cite{E}), was assumed to be 1/2 or 27/40 by Din\'is and Parrondo \cite{DP} and to be 1 by Van den Broeck and Cleuren \cite{VC}.  Here we let $\phi$ range over the interval $(0,1]$.

Behrends \cite{B} introduced a stochastic model that includes our (deterministic) model as a special case, and he proved that the sequence of expected (or average) profits is eventually quasiperiodic under certain assumptions.

This paper develops, in the context of the Din\'is--Parrondo model, techniques for analyzing the asymptotic behavior of a piecewise-linear discrete dynamical system, and may therefore be of interest even to readers unfamiliar with Parrondo's paradox.

In Section \ref{prelim} we formulate a piecewise-linear discrete dynamical system for the capital-dependent Parrondo games played collectively according to the greedy strategy; it is parametrized by $(\rho,\phi)\in(0,1)\times(0,1]$.  In Section \ref{Bforever} we show that, in one region of the parameter space (namely, $\phi\le2/3$), game $B$ is eventually played forever, resulting in an asymptotically fair game and, in terms of the discrete dynamical system, a globally asymptotically stable equilibrium.  In Section \ref{periodic} we show that, in the remainder of the parameter space (namely, $\phi>2/3$), there is an initial state that yields a periodic pattern of games, resulting in an asymptotically winning game and, in terms of the discrete dynamical system, a limit cycle.  In fact, in a very small region of the parameter space, there are two limit cycles.  Section \ref{stability} attempts to show that (again assuming  $\phi>2/3$), where there is a unique limit cycle it is asymptotically stable (in fact, it is globally asymptotically stable unless there is an unstable equilibrium).  The proofs of the assertions in Section \ref{stability} are incomplete, so some of our findings are stated as conjectures.  In Section \ref{profit} we confirm the assertions just made concerning the asymptotic profitability of the greedy strategy.

For example, if $(\rho,\phi)=(1/3,1/2)$ (and $\eps=0$), our results show that there is a globally asymptotically stable equilibrium with game $B$ eventually played forever.  This is contrary to a computational result of Din\'is and Parrondo \cite{DP}, who found the pattern $[1,40]$ in this case (i.e., $AB^{40}AB^{40}AB^{40}\cdots$).  The anomaly is likely attributable to roundoff error;  64-bit arithmetic (C++) is insufficient here.

Let us introduce some additional notation.  As defined in the second paragraph above, the game pattern $[1,2]$ stands for $ABBABBABB\cdots$.  It will be useful to have a concise notation for game sequences that are eventually periodic.  We write $ABABBABBABB\cdots$, for example, as $AB\overline{ABB}$, just as one would write the binary expansion of the fraction 5/14 as $0.01\overline{011}$.

With this notation, we can describe the limit cycles that appear in terms of the game patterns.  They are of two types, either one of
$\overline{AB^2}$, $\overline{AB^4}$, $\overline{AB^6}$, \dots,
or one of
$\overline{AB^4AB^2}$, $\overline{AB^6AB^4}$, $\overline{AB^8AB^6}$, \dots.
For further simplicity, we will also denote these patterns by $[1,2]$, $[1,4]$, $[1,6]$, \dots, and by $[1,4,1,2]$, $[1,6,1,4]$, $[1,8,1,6]$, \dots.

\section{Preliminaries}
\label{prelim}

It is well known that a Markov chain $\{X_n\}_{n\ge0}$ with state space $\{0,1,2\}$ underlies the capital-dependent Parrondo games; here $X_n$ represents the player's capital modulo 3 after $n$ games.  When playing game $B$ it evolves according to the one-step transition matrix
\begin{equation}
\label{PBcirc}
\setlength{\arraycolsep}{1.5mm}
{\bm P}_B^\circ:=\left(\begin{array}{ccc}
0&p_0&1-p_0\\
1-p_1&0&p_1\\
p_1&1-p_1&0
\end{array}\right),
\end{equation}
where $p_0:=\rho^2/(1+\rho^2)$ and $p_1:=1/(1+\rho)$ with $0<\rho<1$. When playing game $A$ it evolves according to the one-step transition matrix $\bm P_A^\circ$ defined by the matrix in (\ref{PBcirc}) with $\rho=1$ (i.e., with $p_0=p_1:=1/2$).  The unique stationary distribution $\bm\pi:=(\pi_0,\pi_1,\pi_2)$ of $\bm P_B^\circ$ is given by
$$
\pi_0 = {1+\rho^2\over2(1+\rho+\rho^2)}, \quad \pi_1 = {\rho(1+\rho)\over2(1+\rho+\rho^2)}, \quad \pi_2 ={1+\rho\over2(1+\rho+\rho^2)},
$$
while that of $\bm P_A^\circ$ is $(1/3,1/3,1/3)$.

Given $0<\phi\le1$, consider a large number $N$ of players, of whom $\phi N$ are selected at random.  Everyone in the sample of size $\phi N$ independently plays game $A$ or everyone independently plays game $B$, the choice determined by the strategy. When playing game $B$, each player uses his own capital to determine which coin to toss. Let $x_0(n)$ be the fraction of the players whose capital is divisible by 3 after $n$ turns.  If the players in the sample collectively choose and individually play game $B$, then the expected average profit, conditioned on $x_0(n)$, is equal to
\begin{eqnarray*}
x_0(n) (2 p_0-1) + [1-x_0(n)](2 p_1 -1)
= - x_0(n) {1-\rho^2\over1+\rho^2} + [1-x_0(n)]{1-\rho\over1+\rho},
\end{eqnarray*}
which is nonpositive if and only if $x_0(n)\ge\pi_0$.
Game $A$ always has expected average profit equal to 0. So the strategy of maximizing the expected average profit at each turn can be summarized by the rules ``play game $A$ if $x_0(n) \geq \pi_0$'' and ``play game $B$ if $x_0(n)<\pi_0$.''
In particular, if both games have expected average profit equal to 0, then game $A$ is played.

We investigate the mean-field limit as $N\to\infty$, in which case the model is deterministic.  (We will not try to justify this; the preceding paragraph was included mainly for motivation.)  Let $x_i$ represent the fraction of players whose capital is congruent to $i$ (mod 3) for $i=0,1,2$. Then $x_0+x_1+x_2=1$.  Thus, in the state space defined by
$$
\Delta := \{ (x_0, x_1, x_2) : x_0 \geq 0,\, x_1 \geq 0,\, x_2 \geq 0,\, x_0+x_1+x_2=1 \}
$$
we have a discrete dynamical system given by
$$
(x_0(n+1),x_1(n+1),x_2(n+1))=F(x_0(n),x_1(n),x_2(n)),\qquad n\ge0,
$$
where
$$
F(x_0,x_1,x_2):=\begin{cases}(x_0,x_1,x_2)\bm P_A&\text{if $x_0\ge\pi_0$,}\\
(x_0,x_1,x_2)\bm P_B&\text{if $x_0<\pi_0$,}\end{cases}
$$
and ${\bm P}_A:=(1-\phi)\bm I+\phi\bm P_A^\circ$ and ${\bm P}_B:=(1-\phi)\bm I+\phi\bm P_B^\circ$.  Clearly, the function $F$ is piecewise linear but discontinuous.

Let us define the projection $p$ of $\Delta$ onto $\{(x_0,x_1): x_0\ge0,\,x_1\ge0,\,x_0+x_1\le1\}$ by $p(x_0,x_1,x_2):=(x_0,x_1)$.  It is a one-to-one transformation.  While the trajectory belongs to $\Delta$, it will often be convenient to regard it as belonging to $p(\Delta)$, a subset of the plane.  In fact, we could redefine $F$ in terms of $2\times2$ matrices, but this does not appear to simplify matters.

To study the asymptotic behavior of the system, we will need the spectral representation for matrices $\bm P_A$ and $\bm P_B$.
We note that $\bm \pi$ is also the unique stationary distribution of $\bm P_B$.  The nonunit eigenvalues of ${\bm P}_B$ are given by
$e_1:= 1- \phi + \phi e_1^\circ$ and $e_2:= 1- \phi + \phi e_2^\circ$,
where, with $S:=\sqrt{(1+\rho^2)(1+4\rho+\rho^2)}$,
$$
e_1^\circ :=-{1\over2}+{(1-\rho)S\over2(1+\rho)(1+\rho^2)}\quad \mbox{and} \quad e_2^\circ:=-{1\over2}-{(1-\rho)S\over2(1+\rho)(1+\rho^2)}
$$
are the nonunit eigenvalues of $\bm P_B^\circ$.
Observe that $e_1=0$ if $\phi=\phi_2$, and $e_2=0$ if $\phi=\phi_1$, where
\begin{eqnarray}
\nonumber
\phi_1:={1\over1-e_2^\circ}&=& {2(1+\rho)(1+\rho^2)\over3(1+\rho)(1+\rho^2) + (1-\rho)S}, \\
\label{phi2}
\phi_2:={1\over1-e_1^\circ}&=& {2(1+\rho)(1+\rho^2)\over3(1+\rho)(1+\rho^2)- (1-\rho)S}.
\end{eqnarray}
Since $0>e_1^\circ>-1/2>e_2^\circ>-1$ for $0<\rho<1$, we have $1/2<\phi_1<2/3<\phi_2<1$ for $0<\rho<1$.
Because the nonunit eigenvalues $e_1$ and $e_2$ will play an important role in what follows, we indicate their dependence on $\phi$ in Table \ref{eigenvalues}.

\begin{table}
\begin{center}
\caption{\label{eigenvalues}Dependence on $\phi$ of the nonunit eigenvalues $e_1$ and $e_2$ {of $\bm P_B$.}}\medskip
{\small \begin{tabular}{ccccc}\hline
$0<\phi<\phi_1$   & $e_1>0$ & $e_2>0$ & $|e_1|>|e_2|$ \\
$\phi=\phi_1$     & $e_1>0$ & $e_2=0$ & $|e_1|>|e_2|$ \\
$\phi_1<\phi<2/3$ & $e_1>0$ & $e_2<0$ & $|e_1|>|e_2|$ \\
$\phi=2/3$        & $e_1>0$ & $e_2<0$ & $|e_1|=|e_2|$ \\
$2/3<\phi<\phi_2$ & $e_1>0$ & $e_2<0$ & $|e_1|<|e_2|$ \\
$\phi=\phi_2$     & $e_1=0$ & $e_2<0$ & $|e_1|<|e_2|$ \\
$\phi_2<\phi\le1$ & $e_1<0$ & $e_2<0$ & $|e_1|<|e_2|$ \\
\hline
\end{tabular}}
\end{center}
\end{table}

We define the diagonal matrix $\bm D:={\rm diag}(1,e_1,e_2)$.  Corresponding right eigenvectors (both for $\bm P_B$ and $\bm P_B^\circ$) are
$$
{\bm r}_0:=\left(\begin{array}{c}
1\\
1\\
1
\end{array}\right),\quad
{\bm r}_1:=\left(\begin{array}{c}
(1+\rho)(1-\rho^2-S)\\
2+\rho+2\rho^2+\rho^3+\rho S\\
-(1+2\rho+\rho^2+2\rho^3-S)
\end{array}\right),
$$
$$
{\bm r}_2:=\left(\begin{array}{c}
(1+\rho)(1-\rho^2+S)\\
2+\rho+2\rho^2+\rho^3-\rho S\\
-(1+2\rho+\rho^2+2\rho^3+S)
\end{array}\right).
$$
They are linearly independent,
so we define ${\bm R}:=({\bm r}_0,{\bm r}_1,{\bm r}_2)$ and ${\bm L}:={\bm R}^{-1}$.  The rows of $\bm L$ are left eigenvectors, and the spectral representation gives
\begin{equation}\label{spectral}
{\bm P}_B ^n={\bm R}{\bm D^n}{\bm L},\qquad n\ge0.
\end{equation}

Of course, ${\bm P}_A$ is the special case $\rho=1$ of ${\bm P}_B$, so it follows from (\ref{spectral}) (or a simple induction argument) that
\begin{eqnarray*}
\setlength{\arraycolsep}{1.5mm}
{\bm P}_A^n=\left(\begin{array}{ccc}
1-2d_n&d_n&d_n\\
d_n&1-2d_n&d_n\\
d_n&d_n&1-2d_n
\end{array}\right),\qquad n\ge0,
\end{eqnarray*}
where $d_n:=[1-(1- 3\phi /2)^n]/3$.

\section{A globally asymptotically stable equilibrium for $\phi\le2/3$}
\label{Bforever}

In this section we show that, corresponding to playing game $B$ forever under the greedy strategy, there is a globally asymptotically stable equilibrium when $\phi\le2/3$ and an unstable equilibrium when $2/3<\phi<\phi_2$.  Let
$$
\Delta_A := \{ (x_0, x_1, x_2)\in \Delta : x_0 \geq \pi_0 \}\quad({\rm resp.,\ }
\Delta_B:=\Delta-\Delta_A)
$$
be the set of states at which game $A$ (resp., game $B$) is chosen under the greedy strategy.  In fact, it will be useful to extend this notation considerably.  For example, $\Delta_{ABBA}$ is the subset of $\Delta_A$ such that the first four games played are $ABBA$ (in that order), and $\Delta_{BBA\overline{B}}$ is the subset of $\Delta_B$ such that the complete game sequence is $BBA\overline{B}$.

\begin{proposition}
\label{prop1}
Under the greedy strategy,
game $A$ is chosen for only finitely many consecutive turns, given an initial state in $\Delta_A$.  If $2/3 \leq \phi \leq1$, then, after only one play of game $A$, game $B$ is played.
\end{proposition}

\begin{proof}
After $n$ plays of game $A$, the initial state, say $(x_0, x_1, x_2) \in \Delta_A$, moves to
\begin{equation}\label{effectofP_A^n}
(x_0, x_1, x_2){\bm P}_A^n
=\bigg({1 \over 3},{1 \over 3},{1 \over 3} \bigg) + \bigg(1- {3 \over 2} \phi \bigg)^n \bigg[(x_0,x_1,x_2)-\bigg({1 \over 3},{1 \over 3},{1 \over 3} \bigg) \bigg].
\end{equation}
If $0<\phi <2/3$, the trajectory converges to the limit $(1/3,1/3,1/3)$ along the line segment that connects the initial state $(x_0, x_1, x_2)$
with $(1/3,1/3,1/3)$ as $n\to\infty$. Since $1/3 < \pi_0$ for $0<\rho<1$, we see that, after a finite number of consecutive plays of game $A$, game $B$ is played.  If $2/3 \leq \phi \leq1$, then we have
$$
(x_0, x_1, x_2){\bm P}_A(1,0,0)^\T={1 \over 3} + \bigg(1- {3 \over 2} \phi \bigg) \bigg(x_0 - {1 \over 3} \bigg) \leq {1\over3}<\pi_0,
$$
which means that, after only one play of game $A$ (requiring $x_0\ge\pi_0>1/3$), game $B$ is played.
\end{proof}

\begin{theorem}
\label{thm2}
If $0< \phi \leq 2/3$, wherever the initial state is located, the greedy strategy chooses game $B$ forever except for an initial finite number of turns.  In particular, the discrete dynamical system has a globally asymptotically stable equilibrium, namely $\bm \pi$, the stationary distribution of $\bm P_B$.  If $2/3 < \phi < 1$, game $B$ is chosen forever only when $\phi<\phi_2$, where $\phi_2$ is defined by (\ref{phi2}), and only when the initial state belongs to one of at most four one-dimensional regions, which will be specified below in terms of $\rho$ and $\phi$.
\end{theorem}

\begin{proof}
Recall that $\Delta_{\overline{B}}$ denotes the set of initial states from which game $B$ is played forever.  Once we know that $\Delta_{\overline{B}}$ is eventually reached from any initial state, we have a \textit{linear} discrete dynamical system
$$
(x_0(n+1),x_1(n+1),x_2(n+1))=(x_0(n),x_1(n),x_2(n))\bm P_B,\qquad n\ge0,
$$
or $(x_0(n),x_1(n),x_2(n))=(x_0,x_1,x_2)\bm P_B^n$ for each $n\ge1$,
with $(x_0,x_1,x_2)\in\Delta_{\overline{B}}$.  But since $\bm P_B$ is irreducible and aperiodic, $\bm P_B^n\to\bm\Pi$, where $\bm\Pi$ is the $3\times 3$ matrix with each row equal to $\bm\pi$.  This implies that $(x_0(n),x_1(n),x_2(n))\to\bm\pi$ for all $(x_0,x_1,x_2)\in\Delta_{\overline{B}}$, and this leads to the global asymptotic stability.

Since the trajectory enters $\Delta_B$ eventually by Proposition \ref{prop1} it is enough to consider a trajectory that starts from $\Delta_B$.
 Let $(x_0, x_1,x_2) \in \Delta_B$ and let $\pi_0 (n)$ be the fraction of the players whose capital is divisible by 3 after $n$ plays of game $B$. Then $\pi_0(0)= x_0 <\pi_0$.  From the spectral representation (\ref{spectral}), we have
\begin{eqnarray*}
\pi_0 (n) &=& (x_0, x_1, x_2){\bm P}_B^n(1,0,0)^\T=(x_0, x_1, 1-x_0 -x_1){\bm R}{\bm D}^n{\bm L}(1,0,0)^\T \\
&=&{1  \over  4 (1 + \rho + \rho^2)S} \{ 2 (1+ \rho^2) S \\
&& \quad{} - e_1^n [2 x_0 (1+ \rho + \rho^2)(1+ \rho^2 - S) + 4 x_1 (1+\rho+\rho^2)(1+\rho^2) \\
& & \qquad{} \qquad{} -(1+2\rho+3\rho^2-S)(1+\rho^2)]  \\
&& \quad{} + e_2^n [2 x_0 (1+ \rho + \rho^2)(1+ \rho^2 + S) + 4 x_1 (1+\rho+\rho^2)(1+\rho^2) \\
& & \qquad{} \qquad{} -(1+2\rho+3\rho^2+S)(1+\rho^2)]\},
\end{eqnarray*}
from which it follows that
\begin{equation}\label{pi0-pi0(n)}
\pi_0 - \pi_0(n) = c_1 e_1^n  - c_2 e_2^n,
\end{equation}
where
\begin{eqnarray*}
c_1&:=&
{1+\rho^2\over2S}\bigg[\bigg(1-{S\over1+\rho^2}\bigg)(x_0-\pi_0)+2(x_1-\pi_1)\bigg],\\
c_2&:=&
{1+\rho^2\over2S}\bigg[\bigg(1+{S\over1+\rho^2}\bigg)(x_0-\pi_0)+2(x_1-\pi_1)\bigg].
\end{eqnarray*}
Notice that $c_1 > c_2$ for $x_0 < \pi_0$ and $0<\rho <1$.  Also, $S/(1+\rho^2)>1$ for $0<\rho <1$.

The region $\Delta_{\overline{B}}$ depends on the nonunit eigenvalues of ${\bm P}_B$, so we derive it separately in the seven cases of Table \ref{eigenvalues} as follows.

Case 1. $0<\phi< \phi_1$.
Since $e_1>e_2>0$, if $c_1 \geq 0$ we have $\pi_0(n) < \pi_0$ for all $n\geq 1$.
On the other hand, if $c_1<0$, we can find a positive integer $n$ such that
$(e_2/e_1)^n \leq c_1/c_2 < (e_2/e_1)^{n-1}$,
and we conclude that $\pi_0(n)\ge\pi_0$.  Thus, it follows that
$$
\Delta_{\overline{B}} =\bigg\{ (x_0, x_1, x_2) \in \Delta_B :
\bigg(1- {S \over 1+ \rho^2}\bigg)(x_0-\pi_0)+2(x_1-\pi_1)\ge0\bigg\}.
$$

For $(x_0, x_1, x_2)\in\Delta_B-\Delta_{\overline{B}}$, we have $x_0 < \pi_0$, hence $x_1 <\pi_1$.
Therefore, $(x_0, x_1, x_2)$ moves to state $(y_0, y_1, y_2):=(x_0, x_1, x_2){\bm P}_B$,
which satisfies
$$
y_0 - x_0={\phi[1-(2+\rho)x_0-(1-\rho)x_1] \over 1+\rho}
>{\phi[1-(2+\rho)\pi_0-(1-\rho)\pi_1] \over 1+\rho}
=0
$$
and
\begin{eqnarray*}
y_1 - x_1&=&{\phi[\rho(1+\rho^2)-\rho(1-\rho)x_0-(1+2\rho)(1+\rho^2)x_1 ] \over (1+\rho)(1+\rho^2)} \\
&>&{\phi[\rho(1+\rho^2)-\rho(1-\rho)\pi_0-(1+2\rho)(1+\rho^2)\pi_1 ] \over (1+\rho)(1+\rho^2)}=0.
\end{eqnarray*}
Therefore, in the region $\Delta_B-\Delta_{\overline{B}}$, as long as
game $B$ is played, the trajectory continues to move in the positive $x_0$ and $x_1$ directions and finally reaches the region $\Delta_A$. (It cannot reach $\Delta_{\overline{B}}$ first, and it cannot remain in $\Delta_B-\Delta_{\overline{B}}$ forever, for in either case it would have started in $\Delta_{\overline{B}}$.)
Once it arrives at a certain state $(y_0,y_1,y_2)\in\Delta_A$, it moves along the line segment between $(y_0,y_1,y_2)$ and $(1/3,1/3,1/3)\in\Delta_{\overline{B}}$, proceeding $3\phi/2$ of the way (see (\ref{effectofP_A^n})).  After one or more such jumps from $\Delta_A$, the trajectory will either reach $\Delta_{\overline{B}}$ or return to $\Delta_B-\Delta_{\overline{B}}$. Even if it returns to $\Delta_B-\Delta_{\overline{B}}$, it eventually enters $\Delta_{\overline{B}}$ because it keeps moving in the positive $x_1$ direction while it visits the two regions $\Delta_B-\Delta_{\overline{B}}$ and $\Delta_A$, and from $\Delta_A$ it moves in the positive $x_1$ direction $3\phi/2$ of the way toward 1/3, ensuring that this alternation between $\Delta_B-\Delta_{\overline{B}}$ and $\Delta_A$ cannot go on forever.  See Figure \ref{thm2-fig}.

\begin{figure}
\centering
\includegraphics[width = 340bp]{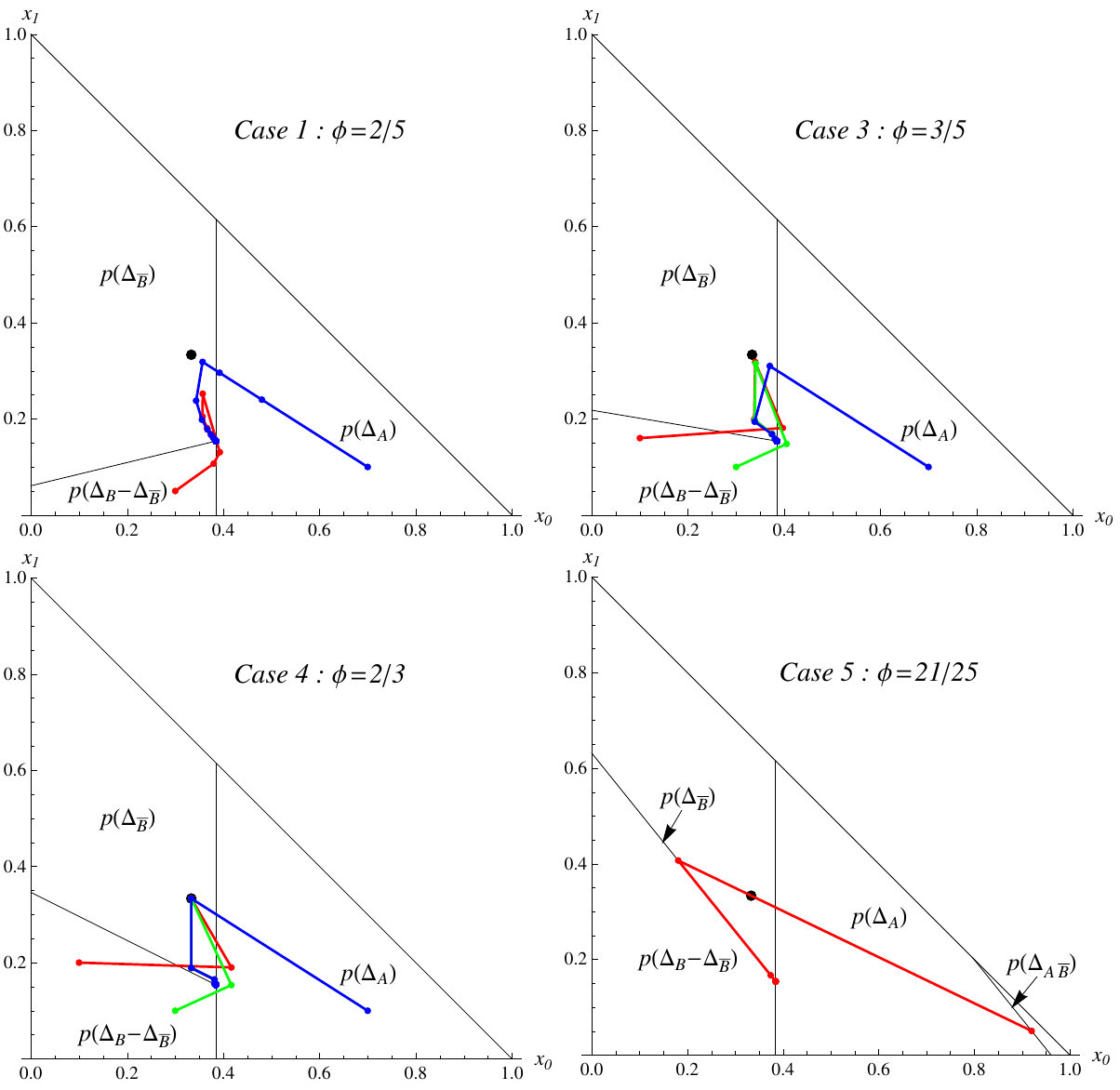}
\caption{\label{thm2-fig}Four of the seven cases of Theorem \ref{thm2} illustrated with $\rho=1/3$.}
\end{figure}

Case 2. $\phi=\phi_1$.
Here $e_1>0=e_2$.
Consequently, $\pi_0(n) < \pi_0$ for all $n\geq 1$ if and only if $c_1 >0$, so that
\begin{eqnarray*}
\Delta_{\overline{B}} = \bigg\{ (x_0, x_1, x_2) \in \Delta_B : \bigg(1- {S \over 1+ \rho^2}\bigg)(x_0-\pi_0)+2(x_1-\pi_1)>0\bigg\}.
\end{eqnarray*}
For $(x_0,x_1,x_2)\in \Delta_B-\Delta_{\overline{B}}$, since $c_1 \leq 0$ we have $\pi_0(1)\ge\pi_0$.  So game $A$ is played after one play of game $B$.
Therefore, when the trajectory starts in the region
$\Delta_B - \Delta_{\overline{B}}$, after a few alternations of game $A$ and game $B$, it enters the region $\Delta_{\overline{B}}$ and stays there forever, just as in Case 1.

Case 3. $\phi_1 < \phi < 2/3$.
In this case, $e_1>0>e_2$ and $|e_1|>|e_2|$.  We claim that a necessary and sufficient condition for $\pi_0(n) < \pi_0$ for every $n\geq 1$ is that $c_1 e_1 - c_2 e_2>0$.  It is clearly necessary (take $n=1$ in (\ref{pi0-pi0(n)})).  For sufficiency, note that $c_1 e_1 - c_2 e_2>0$ implies $c_1>0$.  Thus, we have $c_1 e_1^n-c_2 e_2^n>0$ for all $n\ge1$ if in addition $c_1>c_2\ge0$.  If $c_1>0>c_2$, then for even $n\ge2$, $c_1 e_1^n-c_2 e_2^n>0$, and for odd $n\ge1$,
\begin{equation*}
c_1 e_1^n-c_2 e_2^n=c_1 e_1^n \bigg[1-{c_2\over c_1}\bigg({e_2\over e_1}\bigg)^n\bigg]\ge c_1 e_1^n \bigg( 1- {c_2 e_2 \over c_1 e_1} \bigg)>0,
\end{equation*}
proving the claim.  It follows that
\begin{eqnarray}\label{Delta(B)-case3}
\Delta_{\overline{B}}&=&\bigg\{ (x_0, x_1, x_2) \in \Delta_B : \nonumber\\
&&\qquad\bigg(1+{(3\phi-2)(1+\rho) \over \phi(1-\rho)} \bigg)(x_0 - \pi_0) + 2 (x_1 -\pi_1)>0 \bigg\}.
\end{eqnarray}

Let $(x_0,x_1,x_2) \in \Delta_B -\Delta_{\overline{B}}$.
Since $c_1 e_1 - c_2 e_2 \leq 0$, game $A$ is played after one play of game $B$.
Especially if $x_1<\pi_1$, the trajectory enters $\Delta_{\overline{B}}$ after a few alternations of game $B$ and game $A$
as in Case 1.  If $x_1 \geq \pi_1$, it moves to state
$(y_0, y_1, y_2):=(x_0, x_1, x_2){\bm P}_B$ in $\Delta_A$, which satisfies
\begin{eqnarray*}
y_1-\pi_1
&=& [ 2(1+\rho)(1+\rho^2)(1+\rho+\rho^2)]^{-1} \{[2\phi(1+\rho+\rho^2)-(1+\rho)^2]\rho(1+\rho^2)\\
&&\quad{} - 2\phi \rho(1-\rho^3) x_0 +2[1+\rho-\phi(1+2\rho)](1+\rho^2)(1+\rho+\rho^2)x_1 \} \\
&>& [ 2(1+\rho)(1+\rho^2)(1+\rho+\rho^2)]^{-1} \{[2\phi(1+\rho+\rho^2)-(1+\rho)^2]\rho(1+\rho^2)\\
&&\quad{} - 2\phi \rho(1-\rho^3) \pi_0 +2[1+\rho-\phi(1+2\rho)](1+\rho^2)(1+\rho+\rho^2)\pi_1 \}\\
&=&0,
\end{eqnarray*}
where the inequality uses $\phi <2/3 <(1+\rho)/(1+2 \rho)$ for $0<\rho<1$.
From Proposition \ref{prop1} it follows that the trajectory converges to the limit $(1/3,1/3,1/3)\in \Delta_{\overline{B}}$ along the line segment between
$(y_0, y_1, y_2) \in \Delta_A$ with $y_1 > \pi_1$ and $(1/3,1/3,\break1/3)$. Since this line is above (in $p(\Delta)$) the critical line segment for $\Delta_{\overline{B}}$ in (\ref{Delta(B)-case3}), after a finite number of plays of game $A$ the trajectory directly enters the region $\Delta_{\overline{B}}$.  See Figure \ref{thm2-fig}.

Case 4. $\phi=2/3$.
Since $e_1>0$ and $e_2=-e_1$, for even $n \geq 2$,
$\pi_0 - \pi_0(n) = e_1^n (\pi_0-x_0) >0$,
and for odd $n \geq 1$,
\begin{eqnarray}\label{case4}
\pi_0-\pi_0 (n)=(1+\rho^2)e_1^n[x_0-\pi_0+2(x_1-\pi_1)]/S.
\end{eqnarray}
So we have
\begin{eqnarray*}
\Delta_{\overline{B}} =\{ (x_0, x_1, x_2) \in \Delta_B : x_0-\pi_0+2(x_1-\pi_1)>0\}.
\end{eqnarray*}
If $(x_0, x_1, x_2) \in \Delta_B-\Delta_{\overline{B}}$, since $\pi_0(1)\ge\pi_0$ by (\ref{case4}), game $A$ is played after one play of game $B$.  But after playing game $A$ the trajectory directly moves to $(1/3,1/3,1/3) \in \Delta_{\overline{B}}$ because $\bm P_A$ is the $3\times3$ matrix each of whose entries is 1/3.
Therefore, there are only three possible sequences of games, namely $\overline{B}$, $A\overline{B}$, and $BA\overline{B}$.  See Figure \ref{thm2-fig}.

Case 5. $2/3 < \phi < \phi_2 $.
Here $e_1 >0>e_2$ and $|e_1|<|e_2|$.
If $c_2=0$, then because $c_1>c_2=0$ we have $\pi_0(n) < \pi_0$ for all $n\geq 1$.  On the other hand, if $c_2\ne0$, then $\pi_0(n)\ge\pi_0$ for some $n\ge1$.
This gives
\begin{equation}\label{Delta(B)}
\Delta_{\overline{B}} = \bigg\{ (x_0, x_1, x_2) \in \Delta_B : \bigg(1+{S \over 1+ \rho^2} \bigg)(x_0 - \pi_0) +2(x_1 - \pi_1)=0  \bigg\}.
\end{equation}
We recall from Proposition \ref{prop1} that there is only one play of game $A$ when the trajectory enters $\Delta_A$.
Since $(w_0, w_1, w_2) \in \Delta_A$ moves to
$$
(x_0, x_1, x_2) := (w_0, w_1, w_2){\bm P}_A ={3\over 2}\phi\bigg({1\over3},{1\over3},{1\over3}\bigg)+\bigg(1-{3 \over 2} \phi\bigg)(w_0,w_1,w_2),
$$
any state $(w_0, w_1, w_2) \in \Delta_A$ satisfying
\begin{eqnarray*}
\bigg(1+{S \over 1+ \rho^2} \bigg)\bigg[{1\over2}\phi + \bigg(1-{3 \over 2} \phi\bigg)w_0 - \pi_0 \bigg] +2\bigg[{1\over2}\phi + \bigg(1-{3 \over 2} \phi\bigg)w_1 - \pi_1\bigg]=0
\end{eqnarray*}
moves to $\Delta_{\overline{B}}$ after one play of game $A$.
Defining
\begin{eqnarray*}
f(x)&:=& - {1 \over 2} \bigg(1+ {S \over 1+\rho^2} \bigg) x + { 2( \phi -2 \pi_1 )(1+\rho^2) + (\phi - 2 \pi_0)(1+ \rho^2 + S) \over 2(3 \phi -2)(1+\rho^2) }
\end{eqnarray*}
for $\pi_0 \leq x \leq 1$,
we have
\begin{equation}\label{DeltaA(B)}
\Delta_{A\overline{B}} = \{(w_0, w_1, w_2) \in \Delta_A : w_1 = f(w_0)\}.
\end{equation}
Since the slope $m$ of $y=f(x)$ satisfies $m <-1$, we need both $f(\pi_0)\geq 0$ and $f(1) \leq 0$ to ensure that $\Delta_{A \overline{B}}$ is nonempty. For $0<\rho<1$ and $2/3 < \phi < \phi_2$ it is easy to check that $f(\pi_0)>0$.
Next we have
\begin{eqnarray*}
f(1) ={ (1+\rho)[(1-\rho)(1+\rho^2) + (1+\rho)S] - 2\phi (1+\rho+\rho^2) S  \over 2(3\phi-2)(1+\rho^2)(1+\rho+\rho^2)}.
\end{eqnarray*}
So it follows that if $2/3 < \phi < \phi_3$ with
$$
\phi_3 :={(1+\rho)[(1-\rho)(1+\rho^2) + (1+\rho)S]\over2 (1+\rho+\rho^2) S},
$$
then $f(1) > 0$ and $\Delta_{A \overline{B}}=\varnothing$, whereas if $\phi_3 \leq \phi < \phi_2$, then $\Delta_{A \overline{B}}\ne\varnothing$.  The latter occurs only in regions 5--8 of Figure \ref{region-fig} in Section \ref{stability}.  See Figure \ref{thm2-fig}.

Similar arguments show that
\begin{equation}\label{DeltaBA(B)}
\Delta_{BA\overline{B}}=\{(w_0,w_1,w_2)\in\Delta_B: w_1=g(w_0)\}
\end{equation}
and
\begin{equation}\label{DeltaBBA(B)}
\Delta_{BBA\overline{B}}=\{(w_0,w_1,w_2)\in\Delta_B: w_1=h(w_0)\},
\end{equation}
where
\begin{eqnarray*}
g(x)  &:=& - {1 \over 2} \bigg(1+ {S \over 1+\rho^2} \bigg) x + {g_1 (\rho, \phi)\over (3\phi-2)g_2(\rho, \phi)},
\end{eqnarray*}
$g_1(\rho, \phi) := \phi (3 \phi - 2) [(1+2\rho)(1+\rho^2) + S] -(1 + \rho) [ 2( \phi -2 \pi_1 )(1+\rho^2) + (\phi - 2 \pi_0)(1+ \rho^2 + S)]$,
$g_2 (\rho, \phi) := (3\phi-2)(1+ \rho)(1+\rho^2)  + \phi(1-\rho) S$,
\begin{eqnarray*}
h(x)  &:=& - {1 \over 2} \bigg(1+ {S \over 1+\rho^2} \bigg) x + {h_1 (\rho, \phi) \over (3\phi-2)h_2 (\rho, \phi)},
\end{eqnarray*}
$h_1(\rho, \phi) := \phi (3 \phi - 2) \{2 (1 + \rho) [(1+2\rho)(1+\rho^2) + S]
- \phi [2 + 6 \rho + 3 \rho^2 + 4 \rho^3 + 3 \rho^4 + (2 + 2 \rho - \rho^2) S]\}-  (1 + \rho)^2[2 (\phi - 2 \pi_1) (1 + \rho^2) + (\phi - 2 \pi_0) (1 + \rho^2 + S)]$, and $h_2(\rho, \phi) := -2 (1 + \rho)^2 (1 + \rho^2) + 6 \phi (1 + \rho)^2 (1 + \rho^2) - \phi^2 (5 + 10 \rho + 6 \rho^2 + 10 \rho^3 + 5 \rho^4) -
 \phi (3 \phi - 2) (1 - \rho^2) S$.  Further, (\ref{DeltaBA(B)}) is nonempty in a region slightly smaller than the union of regions 7 and 8 of Figure \ref{region-fig} in Section \ref{stability}, while (\ref{DeltaBBA(B)}) is nonempty only in a very small subset of region 8 of that figure.

Finally, it can be shown that $\Delta_{ABA\overline{B}}=\varnothing$, $\Delta_{ABBA\overline{B}}=\varnothing$, and $\Delta_{BBBA\overline{B}}=\varnothing$, implying that there are no other ways in which game $B$ is played forever.  In summary, in the case of $2/3 <\phi<\phi_2$, only when the initial state belongs to (\ref{Delta(B)}), (\ref{DeltaA(B)}), (\ref{DeltaBA(B)}), or (\ref{DeltaBBA(B)}), some of which may be empty, is game $B$ played forever.

Cases 6 and 7.  $\phi_2\le \phi\le1$.  Here $e_2<e_1\le0$.  Clearly, $\Delta_{\overline{B}}=\varnothing$ in these cases.
\end{proof}

\section{A limit cycle for $\phi>2/3$}
\label{periodic}

In this section we show that, whenever $\phi>2/3$, there is at least one limit cycle, at least when the discrete dynamical system starts from a certain initial state, which will be specified.  The game patterns that occur include $[1,n]$, denoting 1 play of game $A$ followed by $n$ plays of game $B$, for even $n\ge2$, and $[1,n,1,n-2]$, denoting 1 play of game $A$, $n$ plays of game $B$, 1 play of game $A$, and $n-2$ plays of game $B$, for even $n\ge4$. The value of $n$ depends on $\rho$ and $\phi$.

We will need three lemmas to prepare for the next
theorem.  The proofs are trivial and therefore omitted.

\begin{lemma}\label{Lemma1}
If $1<a<b$ and $c>0$, then $(a^n+c)/(b^n+c)$ is decreasing in $n\ge1$.
\end{lemma}

\begin{lemma}\label{Lemma2}
If $0<a<b<1$, $a+b\le1$, and $c>1$, then $(c-a^n)/(c-b^n)$ is
decreasing in $n\ge1$.
\end{lemma}

\begin{lemma}\label{Lemma3}
If $0<c\le{1\over2}$, the functions
$$
f(x):={x-cx^n\over1-cx^n}\quad{\rm and}\quad g(x):={x-cx^n\over1-cx^{n+1}}
$$
are increasing on $(0,1)$ for each $n\ge1$.
\end{lemma}

We will also need the functions
\begin{eqnarray}\label{E_n}
E_n &:=&\phi (1- \rho)\big\{ e_2^n[2+ e_1^n (3 \phi -2) ][3(1+\rho)(1+\rho^2)-(1-\rho)S] \nonumber\\
&& \qquad\qquad\;\;{} - e_1^n [2+ e_2^n (3\phi -2)][3(1+\rho)(1+\rho^2) +(1-\rho)S ] \big\}/2,\\
E_{n,m} &:=&\phi (1- \rho)\big\{ e_2^m[2+ e_1^n (3 \phi -2) ][3(1+\rho)(1+\rho^2)-(1-\rho)S] \nonumber\\ \label{E_{n,m}}
&& \qquad\qquad\;\;{} - e_1^m [2+ e_2^n (3\phi -2)][3(1+\rho)(1+\rho^2) +(1-\rho)S ] \big\}/2,\\
G_{n,m}&:=& \phi (1- \rho)\big\{ e_2^m [4 -e_1^{2 (n-1)}(3 \phi-2)^2][ 2  - e_2^{n-2} (3 \phi -2) ][3(1+\rho)(1+\rho^2)\nonumber\\
&& \qquad\qquad\;\; {}-(1-\rho)S ] - e_1^m [4 -e_2^{2 (n-1)}(3 \phi-2)^2] [ 2 - e_1^{n-2} (3\phi -2) ]\nonumber\\ \label{G_{n,m}}
&&\qquad\qquad\qquad\qquad{}\cdot[3(1+\rho)(1+\rho^2)+(1-\rho)S ] \big\},\\
H_{n,m}&:=&  \phi (1- \rho)\big\{ e_2^m [4-e_1^{2 (n-1)}(3 \phi-2)^2][2 - e_2^n (3 \phi -2)]
[3(1+\rho)(1+\rho^2)\nonumber\\
&& \qquad\qquad\;\; {}-(1-\rho)S ] - e_1^m [4 -e_2^{2 (n-1)}(3 \phi-2)^2] [ 2- e_1^n (3 \phi -2) ]\nonumber\\ \label{H_{n,m}}
&&\qquad\qquad\qquad\qquad{}\cdot[3(1+\rho)(1+\rho^2) +(1-\rho)S ] \big\}.
\end{eqnarray}
The significance of these functions is explained in the following theorem.

\begin{theorem}\label{thm-limitcycles}
For even $n\ge4$, the implicitly defined curves $G_{n,n-2}=0$, $E_{n-2}=0$, $E_{n,n-2}=0$,  $H_{n,n-2}=0$, and $G_{n+2,n}=0$ are monotonically ordered, from highest to lowest, in $\{(\rho,\phi)\in(0,1)\times (2/3,3/4): \phi<\phi_2\}$.  More precisely,
\begin{eqnarray}
G_{n,n-2}=0&{\rm\ is\ above\ }&E_{n-2}=0,\label{ineq1}\\
E_{n-2}=0&{\rm\ is\ above\ }&E_{n,n-2}=0,\label{ineq2}\\
E_{n,n-2}=0&{\rm\ is\ above\ }&H_{n,n-2}=0,\label{ineq3}\\
H_{n,n-2}=0&{\rm\ is\ above\ }&G_{n+2,n}=0,\label{ineq4}
\end{eqnarray}
for $n=4,6,8,\ldots$.  Furthermore, the functions defining the curves are positive above, and negative below, the curves.

For even $n\ge2$ and $(\rho,\phi)\in(0,1)\times(2/3,1]$, the greedy strategy leads to the periodic pattern $[1,n]$ starting from the corresponding stationary distribution as the initial state if and only if $E_{n,n-2}<0$ and $E_n\ge0$.  ($E_{2,0}<0$ is automatically satisfied.)

For even $n\ge4$ and $(\rho,\phi)\in(0,1)\times(2/3,1]$, the greedy strategy leads to the periodic pattern $[1,n,1,n-2]$ starting from the corresponding stationary distribution as the initial state if and only if $G_{n,n-2}<0$ and $H_{n,n-2}\ge0$.
\end{theorem}

\begin{remark}
In particular, if $\phi>2/3$, there is at least one limit cycle.

If $(\rho,\phi)$ belongs to the region between the curves in (\ref{ineq1}), there are two limit cycles, of the forms $[1,n,1,n-2]$ and $[1,n-2]$.  If $(\rho,\phi)$ belongs to the region between the curves in (\ref{ineq2}), there is one limit cycle, of the form $[1,n,1,n-2]$.  If $(\rho,\phi)$ belongs to the region between the curves in (\ref{ineq3}), there are two limit cycles, of the forms $[1,n,1,n-2]$ and $[1,n]$.  If $(\rho,\phi)$ belongs to the region between the curves in (\ref{ineq4}), there is one limit cycle, of the form $[1,n]$.

The regions with two limit cycles are very small.  See Table \ref{regions-rho=1/3} for the case $\rho=1/3$.
\end{remark}

\begin{table}
\caption{\label{regions-rho=1/3}Critical $\phi$-values separating regions at $\rho=1/3$.  Numbers are truncated (not rounded) at 18 decimal places.\medskip}
\begin{center}
{\small \begin{tabular}{ccc}\hline
   form of   & lower-boundary &  $\phi$-value \\
limit cycles & equation &  at $\rho=1/3$ \\
\hline
$[1,2]$              & $G_{4,2}=0$ & $0.688\,066\,413\,565\,052\,628\cdots$ \\
$[1,4,1,2]$, $[1,2]$ & $E_2=0$ & $0.688\,066\,239\,503\,137\,641\cdots$ \\
$[1,4,1,2]$          & $E_{4,2}=0$ & $0.688\,026\,898\,650\,299\,426\cdots$ \\
$[1,4,1,2]$, $[1,4]$ & $H_{4,2}=0$ & $0.688\,026\,881\,018\,074\,821\cdots$ \\
$[1,4]$              & $G_{6,4}=0$ & $0.677\,218\,563\,694\,275\,305\cdots$ \\
$[1,6,1,4]$, $[1,4]$ & $E_4=0$ & $0.677\,218\,563\,614\,298\,209\cdots$ \\
$[1,6,1,4]$          & $E_{6,4}=0$ & $0.677\,217\,953\,395\,292\,912\cdots$ \\
$[1,6,1,4]$, $[1,6]$ & $H_{6,4}=0$ & $0.677\,217\,953\,388\,847\,194\cdots$ \\
$[1,6]$              & $G_{8,6}=0$ & $0.673\,669\,128\,225\,600\,196\cdots$ \\
$\vdots$ & $\vdots$ & $\vdots$ \\
\hline
\end{tabular}}
\end{center}
\end{table}

\begin{proof}
Let $\bm \pi_{[1,n]}$ be the stationary distribution of $ {\bm P}_A {\bm P}_B^n$. Then we have
\begin{eqnarray}\label{E_n/D_n}
(\bm \pi_{[1,n]} - \bm \pi)(1,0,0)^\T = E_n / D_n ,
\end{eqnarray}
where $E_n$ is as in (\ref{E_n}) and
\begin{equation}\label{Dn}
D_n := 2[2+ e_1^n (3 \phi-2)][2+ e_2^n (3 \phi -2)](1+\rho+\rho^2)S.
\end{equation}
Now $D_n$ is positive for all positive integers $n$ because $0<3\phi-2\le1$ and $|e_1|<|e_2|<1$.  Noting that $e_1+e_2=-(3\phi-2)$, $e_2-e_1=-\phi(1-\rho)S/[(1+\rho)(1+\rho^2)]$, and $e_1e_2= 1- 3 \phi + 2 \phi^2 (1+\rho+\rho^2)^2/[(1+\rho)^2(1+\rho^2)]$, we have
$$
E_2 = { \phi^2 (1-\rho)^2 g(\rho, \phi) S \over (1+\rho)^4 (1+ \rho^2)^2},
$$
where $g(\rho, \phi) :=  g_1(\phi)(1 + 4 \rho +4 \rho^7 + \rho^8)+ g_2(\phi) \rho^2 (1 + \rho^4)+ g_3(\phi)\rho^3 (1 + \rho^2)+ g_4(\phi) \rho^4$ with $g_1(\phi):=-15 + 48 \phi - 63 \phi^2 + 44 \phi^3 - 12 \phi^4$,
$g_2(\phi):=-120 + 396 \phi - 540 \phi^2 + 404 \phi^3 -120 \phi^4$,
$g_3(\phi):=-180 + 600 \phi - 828 \phi^2 +632 \phi^3 - 192 \phi^4$, and
$g_4(\phi):=-210 + 696 \phi -954 \phi^2 + 728 \phi^3 - 228 \phi^4$.
Since the functions $g_1$, $g_2$, $g_3$, and $g_4$ are increasing on $(0,1]$, we have
$g(\rho,\phi)\geq g(\rho,\phi_2)$ for $\phi_2  \leq \phi \leq 1$ and
\begin{eqnarray*}
g(\rho,\phi_2) ={16 (1-\rho) (1+\rho)^3(1+\rho^2)^3 g_5(\rho)S \over [3(1+\rho)(1+\rho^2)-(1-\rho)S]^4},
\end{eqnarray*}
where $g_5(\rho):=17 + 68 \rho + 80 \rho^2 + 92 \rho^3 + 134 \rho^4 + 92 \rho^5 + 80 \rho^6+ 68 \rho^7 + 17 \rho^8 -3(1-\rho^2) (5 + 10 \rho +6 \rho^2 + 10 \rho^3 + 5 \rho^4)S$.
Using
\begin{eqnarray*}
&&[17 + 68 \rho + 80 \rho^2 + 92 \rho^3 + 134 \rho^4 + 92 \rho^5 + 80 \rho^6+ 68 \rho^7 + 17 \rho^8]^2\\
   &&\quad{}\quad{} -[3(1-\rho^2)(5 + 10 \rho +6 \rho^2 + 10 \rho^3 + 5 \rho^4)S]^2=64(1+\rho+\rho^2)^8>0,
\end{eqnarray*}
we have $g_5(\rho)>0$ and hence $g(\rho,\phi_2)>0$ for $0<\rho<1$.
It follows that if $\phi_2  \leq \phi \leq 1$, then $E_2 >0$, which implies that $\bm \pi_{[1,2]} \in \Delta_A$.

More generally, a similar argument shows that
\begin{eqnarray}
\label{EnmDn}
(\bm \pi_{[1,n]}{\bm P}_A {\bm P}_B^m - \bm \pi)(1,0,0)^\T = E_{n,m} / D_n ,
\end{eqnarray}
where $E_{n,m}$ and $D_n$ are as in (\ref{E_{n,m}}) and (\ref{Dn}).  Notice that $E_{n,n}=E_n$.
First, $E_{2,0}<0$ by Proposition \ref{prop1}.  Next,
\begin{eqnarray*}
E_{2,1} = {\phi^2 (1-\rho)^2 S\over(1+\rho)^2(1+\rho^2)} [h(\phi)(1+2\rho+2\rho^3+\rho^4)-6(1-\phi)(5-9\phi +9 \phi^2)\rho^2],
\end{eqnarray*}
where $h(\phi):=-15+40\phi-45\phi^2+18\phi^3$.  Since $h$ is increasing, its maximum on $(0,1]$ is $h(1)=-2$, so we have $E_{2,1}<0$, which proves that $\bm \pi_{[1,2]}{\bm P}_A {\bm P}_B \in \Delta_B$.  We need one more play of game $B$ to return to $\bm\pi_{[1,2]}{\bm P}_A {\bm P}^2_{B}=\bm\pi_{[1,2]}\in\Delta_A$.  Hence if $\phi_2\le\phi\le1$, the trajectory that starts from initial state $\bm \pi_{[1,2]} \in \Delta_A$ follows the periodic pattern $[1,2]$ under the greedy strategy.

Next consider the case of $ 2/3 < \phi < \phi_2 $.
We can rewrite $E_n$ in (\ref{E_n}) as
\begin{eqnarray*}
E_n &=& {\phi (1- \rho)\over2}[3(1+\rho)(1+\rho^2) +(1-\rho)S]e_2^n [2+ e_1^n (3\phi -2)] \bigg( {\phi_1 \over\phi_2 }- F_n\bigg),
\end{eqnarray*}
where $F_n :=(e_1/e_2)^n [2+e_2^n(3\phi-2)]/[2+e_1^n(3\phi-2)]$.
In this case, since $e_2<0<e_1$ with $|e_2|>|e_1|$, we find, for all odd $n\ge1$, that $F_n<0$ and therefore $E_n <0$.
Moreover, $F_n>0$ for all even $n\ge2$, and the sequence $\{F_n: n=2,4,6,\ldots\}$ is decreasing to 0 since
$F_n=[(1/e_2)^n+(3\phi-2)/2]/[(1/e_1)^n+(3\phi-2)/2]$
and therefore Lemma \ref{Lemma1} applies.
Thus, we let $s$ denote the smallest even $n\ge2$ such that
$0< F_n \le \phi_1/\phi_2$.  Equivalently, $s:=\min\{n\in\{2,4,6,\ldots\}:E_n\ge0\}$, hence $E_s\ge0$ and (if $s\ge4$) $E_{s-2}<0$.  By (\ref{E_n/D_n}), $\bm \pi_{[1,s]} \in \Delta_A$.

For odd $m<s$, we have $E_{s,m}<0$.  If $E_{s,s-2}<0$,
then $E_{s,m}<0$ for all even $m<s-2$ and since $\bm \pi_{[1,s]} {\bm P}_A{\bm P}_B^s = \bm \pi_{[1,s]}$,
we can conclude that after $s$ plays of game $B$, the trajectory returns to the initial state $\bm \pi_{[1,s]} \in \Delta_A$.
Hence for $(\rho,\phi)$ satisfying $E_{s,s-2}<0$ (with $s$ defined as above), the trajectory that starts from the initial state $\bm \pi_{[1,s]}$ follows the periodic pattern $[1,s]$ under the greedy strategy. (Notice that if $s=2$, then $E_{s,s-2}=E_{2,0}<0$ automatically.)  Moreover, the conclusion fails if $E_n<0$ (implying $\bm\pi_{[1,n]}\in\Delta_B$) or if $E_{n,n-2}\ge0$ (implying $\bm\pi_{[1,n]}\bm P_A\bm P_B^{n-2}\in\Delta_A$).
This proves the assertions in the second paragraph of the theorem.

We next claim that for $(\rho,\phi)$ satisfying $E_{s,s-2} \geq 0$ with $s\geq4$, the greedy strategy leads
to periodic pattern $[1,s,1,s-2]$ if we start from the initial state $\bm \pi_{[1,s,1,s-2]} \in \Delta_A$.  Calculations similar to (\ref{EnmDn}) give
\begin{eqnarray*}
(\bm \pi_{[1,n,1,n-2]}{\bm P}_A {\bm P}_B^m - \bm \pi)(1,0,0)^\T = G_{n,m}/I_n,\\
(\bm \pi_{[1,n,1,n-2]}{\bm P}_A {\bm P}_B^n {\bm P}_A {\bm P}_B^m - \bm \pi)(1,0,0)^\T = H_{n,m}/I_n,
\end{eqnarray*}
where $G_{n,m}$ and $H_{n,m}$ are as in (\ref{G_{n,m}}) and (\ref{H_{n,m}}) and
\begin{eqnarray*}
I_n &:=& 2 [2-e_1^{n-1}(3\phi-2)][2-e_2^{n-1}(3\phi-2)] D_{n-1}.
\end{eqnarray*}
Since $0<3 \phi-2\le1$ and $|e_1|<|e_2|<1$, we have $I_n >0$.
To prove the claim we need to show that
$G_{s,m}<0$ for $0 \leq m \leq s-1$ and $G_{s,s}\ge0$, as well as
$H_{s,m}<0$ for $0 \leq m \leq s-3$ and $H_{s,s-2}\ge0$.
Since $e_2<0<e_1$ with $|e_2|>|e_1|$, it is sufficient to prove that $G_{s,s-2} <0$, $G_{s,s} \geq 0$, $H_{s,s-4} <0$, and $H_{s,s-2} \geq 0$.  From the fact that $E_{s-2}<0$ we have
\begin{eqnarray}\nonumber
\bigg({e_1 \over e_2} \bigg)^{s-2} &>& \bigg({ 2+ e_1^{s-2} (3\phi -2) \over 2+ e_2^{s-2} ( 3 \phi -2)} \bigg) \bigg({\phi_1\over\phi_2 } \bigg) \\ \label{gs-2}
&=& \bigg({ 4- e_1^{2(s-2)} (3\phi -2)^2 \over 4- e_2^{2(s-2)} ( 3 \phi -2)^2} \bigg) \bigg( { 2- e_2^{s-2} (3\phi -2) \over 2- e_1^{s-2}( 3 \phi -2)} \bigg) \bigg({\phi_1 \over\phi_2 } \bigg).
\end{eqnarray}
Here we can show that the sequence $(4- e_1^{2 m}(3\phi-2)^2)/(4-e_2^{2 m}(3 \phi-2)^2)$ is
decreasing in $m$ if we apply Lemma \ref{Lemma2} with
$a=e_1^2$, $b=e_2^2$, and $c=4/(3\phi-2)^2$,
where we recall that $2/3<\phi<\phi_2<1$, so $0<a<b<1$ (see Table \ref{eigenvalues}) and $c>4$.
It remains to show that $a+b\le1$.  Let us write
$e_1=1-\phi+\phi e_1^\circ=1-(3/2)\phi+\phi e^\circ$ and
$e_2=1-\phi+\phi e_2^\circ=1-(3/2)\phi-\phi e^\circ$,
where $e^\circ:=(1-\rho)S/[2(1+\rho)(1+\rho^2)]\in(0,1/2)$.  Then
\begin{eqnarray*}
a+b&=&e_1^2+e_2^2=2\bigg(1-{3\over2}\phi\bigg)^2+2\phi^2(e^\circ)^2<2\bigg(1-{3\over2}\phi\bigg)^2+{\phi^2\over2}\\
&=&2-6\phi+5\phi^2=1+(1-\phi)(1-5\phi)<1
\end{eqnarray*}
since $2/3<\phi<\phi_2<1$.

Hence from (\ref{gs-2}) we have
\begin{eqnarray*}
\bigg({e_1 \over e_2} \bigg)^{s-2}>\bigg({ 4- e_1^{2(s-1)} (3\phi -2)^2 \over 4- e_2^{2(s-1)} ( 3 \phi -2)^2} \bigg)
\bigg( { 2- e_2^{s-2} (3\phi -2) \over 2- e_1^{s-2}( 3 \phi -2)} \bigg) \bigg({\phi_1 \over\phi_2 } \bigg),
\end{eqnarray*}
which is equivalent to $G_{s,s-2}<0$.  By the same reasoning, for even $n\ge4$, $E_{n-2}\le0$ implies $G_{n,n-2}<0$, which yields (\ref{ineq1}).

We also have
\begin{eqnarray}\label{hs-4}
\bigg( {e_1 \over e_2} \bigg) ^{s-4} &=& \bigg( {e_1 \over e_2} \bigg) ^{s-2} \bigg( {e_2 \over e_1} \bigg) ^2 \nonumber\\
&>& \bigg({ 4- e_1^{2(s-1)} (3\phi -2)^2 \over 4- e_2^{2(s-1)} ( 3 \phi -2)^2} \bigg)
\bigg( { 2- e_2^{s-2} (3\phi -2) \over 2- e_1^{s-2}( 3 \phi -2)} \bigg) \bigg({\phi_1 \over\phi_2 } \bigg)\bigg( {e_2 \over e_1} \bigg) ^2 \nonumber\\
&>& \bigg({ 4- e_1^{2(s-1)} (3\phi -2)^2 \over 4- e_2^{2(s-1)} ( 3 \phi -2)^2} \bigg) \bigg( { 2- e_2^s (3\phi -2) \over 2- e_1^s( 3 \phi -2)} \bigg) \bigg({\phi_1 \over\phi_2 } \bigg),
\end{eqnarray}
which is equivalent to $H_{s,s-4}<0$.
Notice that the last inequality in (\ref{hs-4}) uses
\begin{eqnarray}\label{inequality}
{ 2 e_2^2- e_2^s (3\phi -2) \over 2 e_1^2- e_1^s ( 3 \phi -2)}  > { 2- e_2^s (3\phi -2) \over 2- e_1^s( 3 \phi -2)},
\end{eqnarray}
which is equivalent to
$$
{2 e_2^2- e_2^s (3\phi -2) \over2- e_2^s (3\phi -2) }  > {2 e_1^2- e_1^s ( 3 \phi -2)  \over 2- e_1^s( 3 \phi -2)}
$$
for even $s\ge4$.  To confirm the latter inequality, divide both numerators and denominators by 2 and apply Lemma \ref{Lemma3}.

On the other hand, from $E_{s,s-2}\geq 0$ and the argument below (\ref{gs-2}), it follows that
\begin{eqnarray}\label{H-ineq}
\bigg( {e_1 \over e_2} \bigg)^{s-2} &\leq& \bigg({2+e_1^s (3\phi-2)\over 2+ e_2^s (3\phi-2)} \bigg) \bigg({\phi_1\over\phi_2} \bigg) \nonumber\\
&=& \bigg({ 4- e_1^{2s} (3\phi -2)^2 \over 4-e_2^{2s} (3\phi-2)^2} \bigg) \bigg( {2-e_2^s (3\phi-2)\over 2-e_1^s(3\phi-2)} \bigg) \bigg({\phi_1 \over\phi_2} \bigg) \nonumber\\
&<& \bigg({4-e_1^{2(s-1)} (3\phi-2)^2 \over 4-e_2^{2(s-1)} (3\phi-2)^2} \bigg) \bigg({2-e_2^s (3\phi-2) \over 2-e_1^s(3\phi -2)} \bigg) \bigg({\phi_1 \over\phi_2 } \bigg),\qquad
\end{eqnarray}
which is equivalent to $H_{s,s-2}>0$. By the same reasoning, for even $n\ge4$, $E_{n,n-2}\ge0$ implies $H_{n,n-2}>0$, which yields (\ref{ineq3}).

Similarly, we have
\begin{eqnarray}\label{G-ineq}
\bigg( {e_1 \over e_2} \bigg) ^s &=& \bigg( {e_1\over e_2} \bigg) ^{s-2} \bigg( {e_1\over e_2} \bigg) ^2 \nonumber\\
&<&\bigg({4-e_1^{2(s-1)} (3\phi-2)^2 \over 4-e_2^{2(s-1)} (3\phi-2)^2} \bigg)
\bigg({ 2-e_2^s (3\phi -2) \over 2-e_1^s (3\phi-2)} \bigg)
\bigg({\phi_1 \over\phi_2} \bigg)\bigg( {e_1 \over e_2} \bigg) ^2  \nonumber\\
&<& \bigg({ 4-e_1^{2(s-1)} (3\phi-2)^2 \over 4-e_2^{2(s-1)} (3\phi-2)^2} \bigg)
\bigg({2-e_2^{s-2} (3\phi-2) \over 2-e_1^{s-2} (3\phi-2)} \bigg) \bigg({\phi_1 \over\phi_2 } \bigg),
\end{eqnarray}
where the last inequality follows from (\ref{inequality}), and this is equivalent to $G_{s,s}>0$.  This almost proves the assertions in the third paragraph of the theorem, except that we have implicitly assumed that $E_{n,n-2}\ge0$ and $E_{n-2}<0$, both of which are stronger than necessary.  We can weaken the former to $H_{n,n-2}\ge0$, in which case (\ref{H-ineq}) is no longer necessary and (\ref{G-ineq}) follows as before.  We can weaken the latter to $G_{n,n-2}<0,$ in which case (\ref{gs-2}) is no longer necessary and (\ref{hs-4}) follows as before.  (When $E_{n-2}\ge0$ and $G_{n,n-2}<0$, we have $s=n-2$, so we apply the inequalities involving $s$ with $s$ replaced by $n$.)  Finally, the necessity of the inequalities is clear:  If $G_{n,n-2}\ge0$, then $\bm\pi_{[1,n,1,n-2]}\bm P_A\bm P_B^{n-2}\in\Delta_A$, while if $H_{n,n-2}<0$, then $\bm\pi_{[1,n,1,n-2]}=\bm\pi_{[1,n,1,n-2]}\bm P_A\bm P_B^n\bm P_A\bm P_B^{n-2}\in\Delta_B$.

For (\ref{ineq2}), it is enough to show that, for even $n\ge4$, $E_{n,n-2}-E_{n-2}>0$ for all $(\rho,\phi)\in(0,1)\times(2/3,3/4)$.  This will then imply that, for even $n\ge4$, when $E_{n-2}=0$ we have $E_{n,n-2}>0$, which yields (\ref{ineq2}).  Now
\begin{eqnarray*}
E_{n,n-2}-E_{n-2}&=&e_1^{n-2} e_2^{n-2} \phi^2 (3 \phi-2) (1-\rho)^2 S
[6(1+\rho)^2 (1+\rho^2) \\
&&\quad{}- \phi (7 + 14 \rho + 12 \rho^2 + 14 \rho^3 + 7 \rho^4)]
/[(1 + \rho)^2 (1 + \rho^2)],
\end{eqnarray*}
which has the sign of $6(1+\rho)^2 (1+\rho^2) - \phi (7 + 14 \rho + 12 \rho^2 + 14 \rho^3 + 7 \rho^4)$.  The latter is decreasing in $\phi$ and, at $\phi=3/4$, equals $(3/4) (1 + 2 \rho + 4 \rho^2 + 2 \rho^3 + \rho^4)>0$, hence it is positive in $(0,1)\times(2/3,3/4)$.

For (\ref{ineq4}), it is enough to show that, for even $n\ge4$, if $H_{n,n-2}=0$, then $G_{n+2,n}>0$.  Equivalently, it suffices to show that, if
$$
\bigg({e_2 \over e_1} \bigg)^{n-2} \bigg({ 4- e_1^{2n-2} (3\phi -2)^2 \over 4-e_2^{2n-2} ( 3 \phi -2)^2} \bigg) \bigg({ 2- e_2^n (3\phi -2) \over 2-e_1^n ( 3 \phi -2)} \bigg)  \bigg({\phi_1\over\phi_2 } \bigg)=1,
$$
then
$$
\bigg({e_2 \over e_1} \bigg)^n \bigg({ 4- e_1^{2n+2} (3\phi -2)^2 \over 4-e_2^{2n+2} ( 3 \phi -2)^2} \bigg) \bigg({ 2- e_2^n (3\phi -2) \over 2-e_1^n ( 3 \phi -2)} \bigg)  \bigg({\phi_1\over\phi_2 } \bigg)>1.
$$
For this we need only show that, for the last two products of fractions, the ratio of the second to the first is greater than 1.  It is in fact equal to
$$
\bigg({ 4- e_1^{2n+2} (3\phi -2)^2 \over 4-e_2^{2n+2} ( 3 \phi -2)^2} \bigg) \bigg({ 4e_2^2 - e_2^{2n} (3\phi -2)^2 \over 4e_1^2 -e_1^{2n} ( 3 \phi -2)^2} \bigg).
$$
This is greater than 1 if and only if
$$
{ 4e_1^2 -e_1^{2n} ( 3 \phi -2)^2 \over 4- e_1^{2n+2} (3\phi -2)^2 } <  { 4e_2^2 - e_2^{2n} (3\phi -2)^2  \over 4-e_2^{2n+2} ( 3 \phi -2)^2}.
$$
Since $2/3<\phi<\phi_2$, we have $e_1^2<e_2^2$ by Table 1.  We divide both numerators and denominators by 4 and apply Lemma \ref{Lemma3}.

Finally, it remains to show that the functions defining the curves in (\ref{ineq1})--(\ref{ineq4}) are positive above, and negative below, the curves.
Consider $G_{n,n-2}$ for even $n\ge4$.  Notice that $G_{n,n-2}$ is positive, 0, or negative according to whether
\begin{equation}\label{product}
\bigg({e_2\over e_1}\bigg)^{n-2}\bigg({4-e_1^{2 (n-1)}(3 \phi-2)^2\over4 -e_2^{2 (n-1)}(3 \phi-2)^2}\bigg) \bigg({ 2  - e_2^{n-2} (3 \phi -2)\over 2 - e_1^{n-2} (3\phi -2)}\bigg) \bigg({\phi_1\over\phi_2 } \bigg)
\end{equation}
is $>1$, $=1$, or $<1$.  It is therefore enough to show that (\ref{product}) is increasing in $\phi\in(2/3,\phi_2)$ for each $\rho\in(0,1)$.  This follows by showing that the product of the first and third factors is increasing and the second factor alone is increasing (the fourth factor is constant).  The other functions, $E_{n-2}$, $E_{n,n-2}$, and $H_{n,n-2}$, are treated similarly.
\end{proof}

We will later need the following consequence of the proof.

\begin{proposition}\label{propTh7}
Using the notation (\ref{E_n})--(\ref{H_{n,m}}), if $E_{n,n-2}<0$ and $E_n\ge0$ for some even $n\ge4$, then $E_{n,m}<0$ for $m=0,1,\ldots,n-1$.  If $G_{n,n-2}<0$ and $H_{n,n-2}\ge0$ for some even $n\ge4$, then $G_{n,m}<0$ for $m=0,1,\ldots,n-1$ and $H_{n,m}<0$ for $m=0,1,\ldots,n-3$.
\end{proposition}

\section{Asymptotic stability of limit cycles}
\label{stability}

In this section, we consider the case in which the discrete dynamical system starts from an arbitrary initial state. We investigate its asymptotic behavior for parameters $(\rho,\phi)$ belonging to $(0,1)\times(2/3,1]$.

Let $(x_0,x_1,x_2)\in\Delta_A$ be the initial state.  (If the initial state is in $\Delta_B$, the trajectory will enter $\Delta_A$ eventually by Theorem \ref{thm2}, with an exception as noted in that theorem.)  Define
$(y_0,y_1,y_2):=(x_0,x_1,x_2)\bm P_A$,
$(z_0,z_1,z_2):=(y_0,y_1,y_2)\bm P_B$, and
$(w_0,w_1,w_2):=(z_0,z_1,z_2)\bm P_B$.
We know by Proposition \ref{prop1} that $(y_0,y_1,y_2)\in\Delta_B$.  In order that $(z_0,z_1,z_2)\in\Delta_B$, we need $z_0<\pi_0$.  Now
$z_0-\pi_0=\alpha_1 x_0+\beta_1 x_1-\gamma_1$, where
\begin{eqnarray*}
\alpha_1&:=&(3 \phi-2) [-(1 + \rho) + \phi (2 + \rho)]/[2 (1 + \rho)]>0,\\
\beta_1&:=&\phi (3 \phi-2) (1 - \rho)/[2 (1 + \rho)]>0,\\
\gamma_1&:=&{(1 + \rho) (1 + \rho^2)  -
  \phi (3 + \rho) (1 + \rho + \rho^2) + 3 \phi^2 (1 + \rho + \rho^2)\over 2 (1 + \rho) (1 + \rho + \rho^2)}>0,
\end{eqnarray*}
and so $(z_0,z_1,z_2)\in\Delta_B$ is equivalent to
$\alpha_1 x_0+\beta_1 x_1<\gamma_1$.
Since
$$
{\gamma_1\over\alpha_1}-\pi_0={\phi [1 - 3(1-\phi)\rho] (1-\rho^2)\over
 2 (3 \phi-2) (1 + \rho + \rho^2) [-(1 + \rho) + \phi (2 + \rho)]}>0,
$$
the region $\{(x_0,x_1,x_2)\in\Delta_A: \alpha_1 x_0+\beta_1 x_1<\gamma_1\}$ is nonempty.

Next, in order that $(w_0,w_1,w_2)\in\Delta_A$ as well, we need $w_0\ge \pi_0$.  Now $w_0-\pi_0=-\alpha_2 x_0-\beta_2 x_1+\gamma_2$, where
\begin{eqnarray*}
\alpha_2&:=&(3\phi-2)[(1+\rho)^2(1+\rho^2)-2\phi(1+\rho)(2+\rho)(1+\rho^2)\\
&&\quad{}+\phi^2 (4+5\rho+3\rho^2+5\rho^3+\rho^4)]/[2(1+\rho)^2(1+\rho^2)]>0,\\
\beta_2&:=&(3 \phi-2)^2 \phi (1 - \rho)/[2 (1 + \rho)]>0,\\
\gamma_2&:=&[-(1 + \rho)^2(1 + \rho^2)^2  + \phi(1 + \rho) (5 + \rho) (1 + \rho^2) (1 + \rho + \rho^2)\\
&&\quad{}- 2 \phi^2 (1 + \rho^2) (5 + 5 \rho - \rho^2) (1 + \rho + \rho^2)\\
&&\quad{}+\phi^3 (1 + \rho + \rho^2) (7 +5 \rho + 3 \rho^2 + 5 \rho^3 - 2 \rho^4)] \\
&&\;/[2 (1 + \rho)^2 (1 + \rho^2) (1 + \rho + \rho^2)],
\end{eqnarray*}
and so $(w_0,w_1,w_2)\in\Delta_A$ is equivalent to
$\alpha_2 x_0+\beta_2 x_1\le\gamma_2$.
However, $\gamma_2$ is not necessarily positive.  More importantly, the set of $(\rho,\phi)\in(0,1)\times(2/3,1]$ for which
\begin{eqnarray*}
{\gamma_2\over\alpha_2}-\pi_0&=&(1 - \rho) \phi [-(1 + \rho)^2 (1 - 5 \rho) (1 + \rho^2)-12\phi\rho (1 + \rho)^2 (1 + \rho^2)\\
\noalign{\vglue-2mm}
&&\qquad\qquad\quad{} + \phi^2(2 + 11 \rho + 20 \rho^2 + 16 \rho^3 + 16 \rho^4 + 7 \rho^5) ]\\
&&\quad/\{2(3\phi-2)(1+\rho+\rho^2)[(1+\rho)^2(1+\rho^2)\\
&&\qquad{}-2\phi(1+\rho)(2+\rho)(1+\rho^2)+\phi^2(4+5\rho+3\rho^2+5\rho^3+\rho^4)]\}\\
&\ge&0,
\end{eqnarray*}
is the region for which $\{(x_0,x_1,x_2)\in\Delta_A:\alpha_2 x_0+\beta_2 x_1\le\gamma_2\}$ is necessarily nonempty, which includes (but is not equal to) the union of regions 1--12 in Figure \ref{region-fig} below.

Let us assume that $(\rho,\phi)$ belongs to this region.  Then
\begin{equation*}
\Delta_{ABBA}=\{(x_0,x_1,x_2)\in\Delta_A: \alpha_1 x_0+\beta_1 x_1<\gamma_1,\,\alpha_2 x_0+\beta_2 x_1\le\gamma_2\}
\end{equation*}
is nonempty, that is, there exists $(x_0,x_1,x_2)\in\Delta_A$ such that the greedy strategy begins with the game sequence $ABBA$.
If we solve the system of equations
\begin{equation}\label{lines}
\alpha_1 x_0+\beta_1 x_1=\gamma_1\quad{\rm and}\quad
\alpha_2 x_0+\beta_2 x_1=\gamma_2,
\end{equation}
simultaneously, we find that $x_0-\pi_0=-\phi (1 - \rho)^2/[2 (3 \phi-2) (1 + \rho + \rho^2)]<0$, hence the two lines intersect only outside $p(\Delta_A)$ and only one of the inequalities $\alpha_1 x_0+\beta_1 x_1<\gamma_1$ and $\alpha_2 x_0+\beta_2 x_1\le\gamma_2$ is needed to define $p(\Delta_{ABBA})$.  This tells us that $p(\Delta_{ABBA})$ is the intersection of the triangular region
$\{(x_0,x_1): x_0\ge\pi_0,\,x_1\ge0,\,x_0+x_1\le1\}$
with exactly one of the two triangular regions
\begin{eqnarray*}
&&\{(x_0,x_1): x_0\ge\pi_0,\,x_1\ge0,\,\alpha_1 x_0+\beta_1 x_1<\gamma_1\},\\
&&\{(x_0,x_1): x_0\ge\pi_0,\,x_1\ge0,\,\alpha_2 x_0+\beta_2 x_1\le\gamma_2\},
\end{eqnarray*}
hence its closure is convex with at most four extreme points.

Now the lines (\ref{lines}) intersect the vertical line $x_0=\pi_0$ at
$x_1=a_1:=(\gamma_1-\alpha_1\pi_0)/\beta_1>0$ and
$x_1=a_2:=(\gamma_2-\alpha_2\pi_0)/\beta_2\ge0$,
and they intersect the horizontal axis $x_1=0$ at
$x_0=b_1:=\gamma_1/\alpha_1>\pi_0$ and
$x_0=b_2:=\gamma_2/\alpha_2\ge\pi_0$.
The quantities $a_1,a_2,b_1,b_2$, which are functions of $\rho$ and $\phi$, play an important role in what follows.  In particular, they allow us to partition the set of all $(\rho,\phi)\in(0,1)\times(2/3,1]$ such that $G_{4,2}>0$
into 12 regions, as indicated in Figure \ref{region-fig}.

To be more precise, region 1 is $G_{4,2}\ge0$, $a_2<1-\pi_0$, and $\phi<(3/4)(1-\rho)+(2/3)\rho$; region 2 is $a_2\ge1-\pi_0$, $b_2<1$, and $\phi<(3/4)(1-\rho)+(2/3)\rho$; region 3 is $b_2\ge1$; region 4 is $b_2<1$, $\phi<\phi_3$, and $\phi>(3/4)(1-\rho)+(2/3)\rho$; region 5 is $\phi\ge\phi_3$ and $b_1>1$; region 6 is $b_1\le1$, $a_2\ge1-\pi_0$, and $b_2<b_1$; region 7 is $a_2<1-\pi_0$ and $\phi<1-\rho/3$; region 8 is $\phi\ge1-\rho/3$ and $b_2<b_1$; region 9 is $b_2\ge b_1$ and $\phi<1-\rho/3$; region 10 is $\phi\ge1-\rho/3$ and $a_2\ge1-\pi_0$; region 11 is $b_2\ge b_1$, $a_2<1-\pi_0$, and $\alpha_1 f_0+\beta_1 f_1<\gamma_1$, where $(f_0,f_1,f_2):=(1,0,0)\bm P_A\bm P_B$; and region 12 is $\alpha_1 f_0+\beta_1 f_1\ge\gamma_1$.

\begin{figure}
\centering
\includegraphics[width = 340bp]{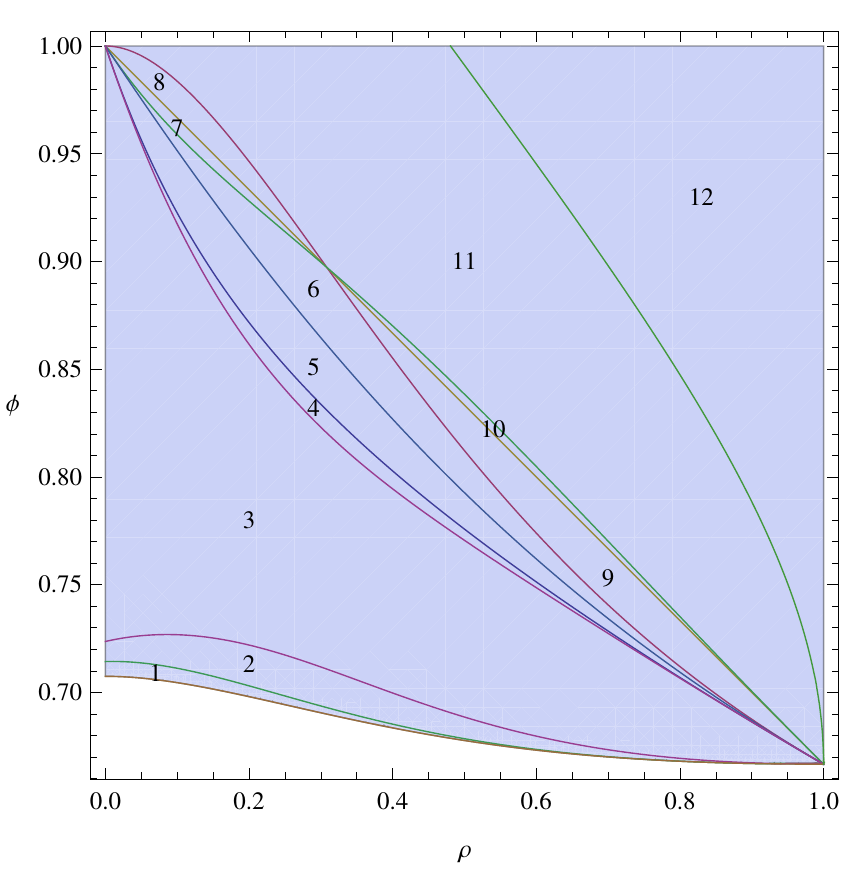}
\caption{\label{region-fig}Partition into 12 regions of the set where $[1,2]$ is (conjectured to be) the form of the unique (asymptotically stable) limit cycle.  If the curves are labeled 1--9 according to their $\phi$-values at $\rho=0.2$ from smallest to largest, curve 1 is $G_{4,2}=0$ ($>$ above), curves 2 and 7 are $a_2=1-\pi_0$ ($>$ between), curves 3 and 4 are $b_2=1$ ($>$ between), curve 5 is $\phi=\phi_3$, curve 6 is $b_1=1$ ($>$ below), curve 8 is $a_1=1-\pi_0$ ($>$ below) or $\phi=1-\rho/3$, and curve 9 is $b_2=b_1$ ($>$ above) or $\phi=\phi_2$.  The curve separating regions 11 and 12 is $\alpha_1 f_0+\beta_1 f_1=\gamma_1$ ($>$ above), where $(f_0,f_1,f_2):=(1,0,0)\bm P_A\bm P_B$.  The regions are defined more precisely in the text.}
\end{figure}

We believe that, whenever $G_{4,2}>0$, there is a unique (asymptotically stable) limit cycle of the form $[1,2]$.  However, we can prove this only in regions 3, 9, 10, and 11.

\begin{theorem}\label{thm-regions3,9--11}
For $(\rho,\phi)$ belonging to region 3 of Figure \ref{region-fig}, there is a unique limit cycle, which is asymptotically stable and of the form $[1,2]$, as well as an unstable equilibrium.  Indeed,
\begin{equation}\label{region3}
\Delta_A=\Delta_{\overline{ABB}},
\end{equation}
and for initial states in $\Delta_B-\Delta_{\overline{B}}$ (see (\ref{Delta(B)})), the trajectory eventually enters $\Delta_A$.

For $(\rho,\phi)$ belonging to region 9, 10, or 11 of Figure \ref{region-fig}, there is a unique limit cycle, which is globally asymptotically stable and of the form $[1,2]$.  Indeed,
\begin{equation}\label{regions9--11}
\Delta_{ABBA}=\Delta_{\overline{ABB}},\qquad \Delta_A-\Delta_{ABBA}=\Delta_{AB\overline{ABB}},
\end{equation}
and for initial states in $\Delta_B$, the trajectory eventually enters $\Delta_A$.
\end{theorem}

\begin{remark}
The theorem includes the case $(\rho,\phi)=(1/3,1)$ studied by Van den Broeck and Cleuren \cite{VC}.
Results for the regions not covered by the theorem will be stated later in the form of a conjecture.
\end{remark}

\begin{proof}
Suppose we can show (with some exceptions in the case of region 3) that the trajectory eventually reaches $\Delta_{\overline{ABB}}$.  Then we have three \textit{linear} discrete dynamical systems
\begin{eqnarray*}
&&(x_0(3m+3),x_1(3m+3),x_2(3m+3))=(x_0(3m),x_1(3m),x_2(3m))\bm P_A\bm P_B\bm P_B,\\
&&(x_0(3m+4),x_1(3m+4),x_2(3m+4))\\
&&\qquad\qquad\qquad\qquad\qquad\quad{}=(x_0(3m+1),x_1(3m+1),x_2(3m+1))\bm P_B\bm P_B\bm P_A,\\
&&(x_0(3m+5),x_1(3m+5),x_2(3m+5))\\
&&\qquad\qquad\qquad\qquad\qquad\quad{}=(x_0(3m+2),x_1(3m+2),x_2(3m+2))\bm P_B\bm P_A\bm P_B,
\end{eqnarray*}
or
\begin{eqnarray*}
(x_0(3m),x_1(3m),x_2(3m))&=&(x_0,x_1,x_2)(\bm P_A\bm P_B\bm P_B)^m,\\
(x_0(3m+1),x_1(3m+1),x_2(3m+1))&=&(x_0,x_1,x_2)(\bm P_A\bm P_B\bm P_B)^m\bm P_A,\\
(x_0(3m+2),x_1(3m+2),x_2(3m+2))&=&(x_0,x_1,x_2)(\bm P_A\bm P_B\bm P_B)^m\bm P_A\bm P_B,
\end{eqnarray*}
with $(x_0,x_1,x_2)\in\Delta_{\overline{ABB}}$.
Now $(\bm P_A\bm P_B\bm P_B)^m\to\bm\Pi_{[1,2]}$, where $\bm\Pi_{[1,2]}$ is the $3\times3$ matrix with each row equal to $\bm\pi_{[1,2]}$, so the limits in the last group of equations are $\bm\pi_{[1,2]}$, $\bm\pi_{[1,2]}\bm P_A$, and $\bm\pi_{[1,2]}\bm P_A\bm P_B$, regardless of the initial state in $\Delta_{\overline{ABB}}$.  This leads to the asymptotic stability.

Recall that $\Delta_{ABBA}$ is the intersection of three triangular regions.  This intersection is itself triangular if
$\min\{a_1,a_2\}\le1-\pi_0$ and $\min\{b_1,b_2\}\le1$
or if
$\min\{a_1,a_2\}\ge1-\pi_0$ and $\min\{b_1,b_2\}\ge1$.
And the intersection is four-sided if
$\min\{a_1,a_2\}>1-\pi_0$ and $\min\{b_1,b_2\}<1$
($\min\{a_1,a_2\}<1-\pi_0$ and $\min\{b_1,b_2\}>1$ is impossible; see Figure \ref{region-fig}).

If for some $(\rho,\phi)$, $\bm P_A\bm P_B^2$ maps the closure of $\Delta_{ABBA}$ into $\Delta_{ABBA}$, then, for this $(\rho,\phi)$, $\Delta_{ABBA}=\Delta_{\overline{ABB}}$.  Now $\bm P_A\bm P_B^2$ is linear, so $\bm P_A\bm P_B^2$ maps the closure of $\Delta_{ABBA}$ into $\Delta_{ABBA}$ if and only if it maps the three or four vertices of the closure of $\Delta_{ABBA}$ into $\Delta_{ABBA}$.  In one important case, this completes the proof.  If $\Delta_{ABBA}=\Delta_A$, then we conclude that (\ref{region3}) holds.  This happens when
$a_1>1-\pi_0$, $a_2\ge1-\pi_0$, $b_1>1$, and $b_2\ge1$.
Actually, the last inequality implies the first three and holds precisely in region 3 of Figure \ref{region-fig}.  Recalling the result of Theorem \ref{thm2}, we conclude that the assertions about region 3 are established.  See Figure \ref{regions3,11-fig}.

\begin{figure}
\centering
\includegraphics[width = 170bp]{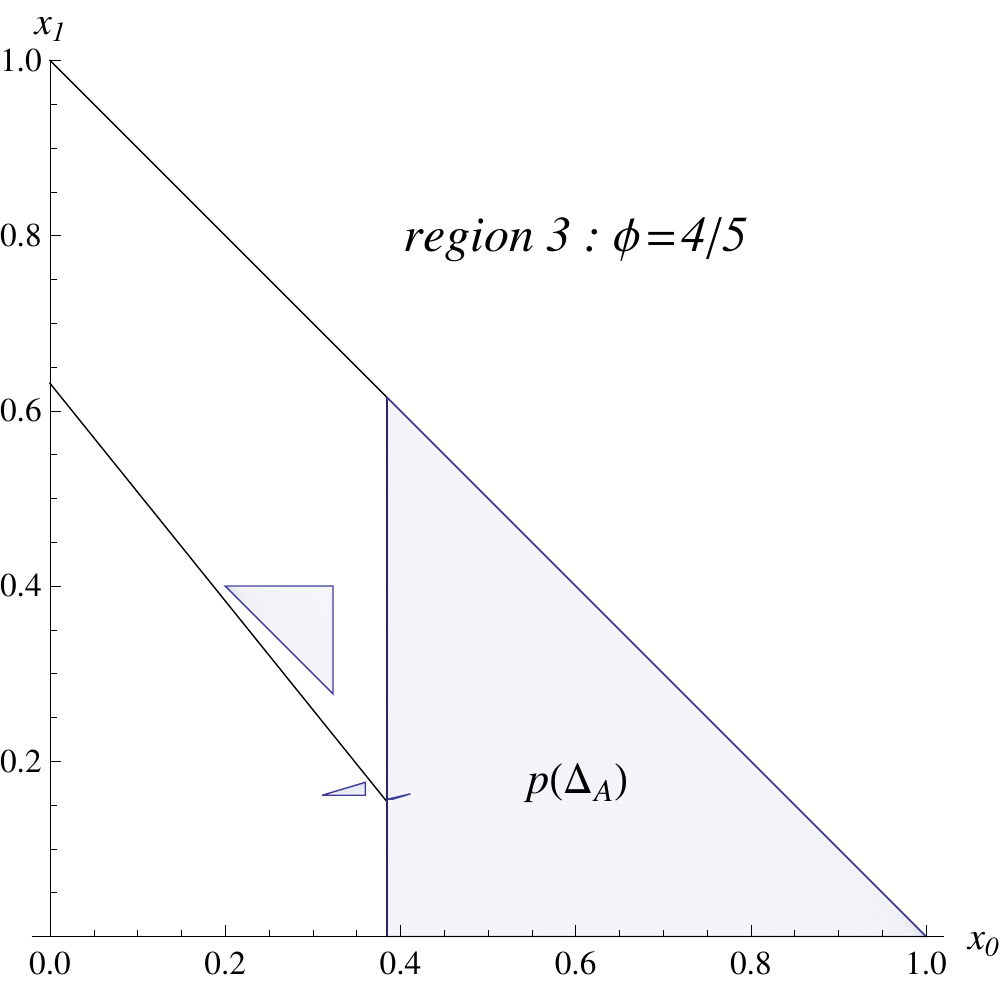}
\includegraphics[width = 170bp]{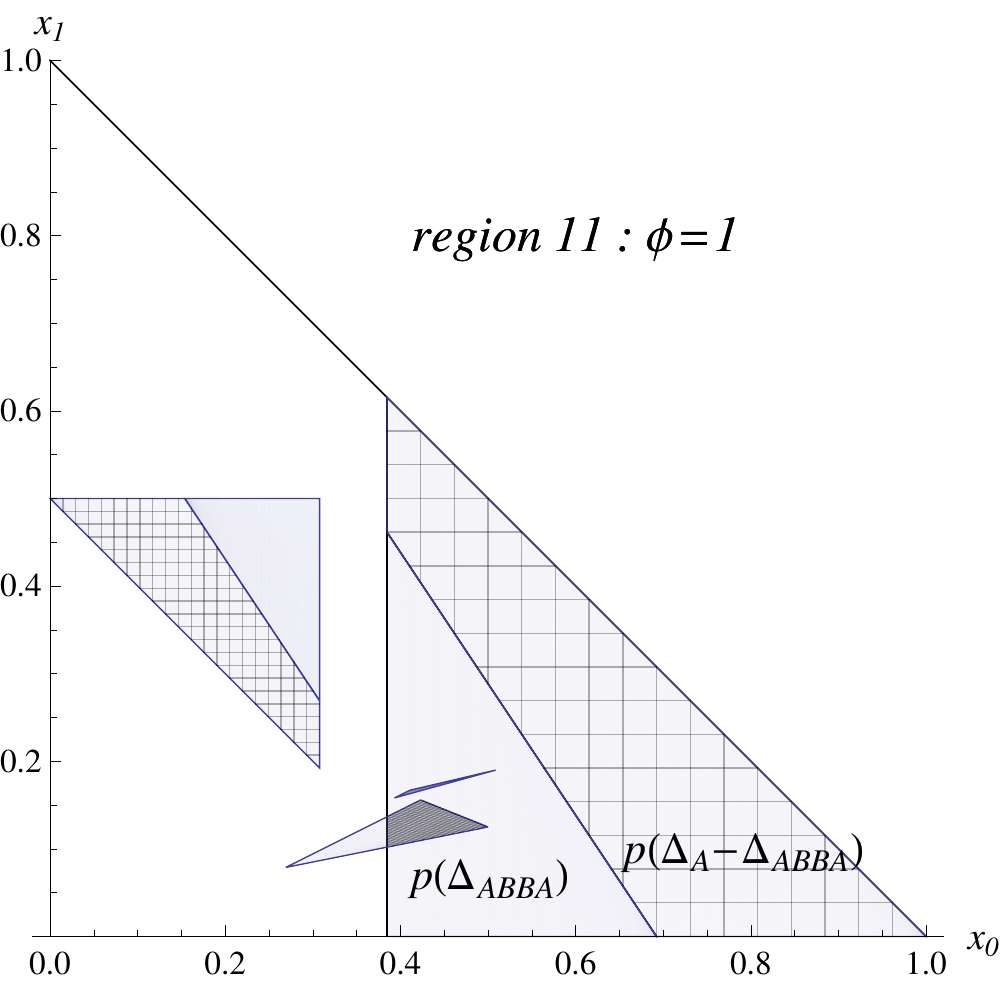}
\caption{\label{regions3,11-fig}Two of the four cases of Theorem \ref{thm-regions3,9--11} illustrated with $\rho=1/3$.}
\end{figure}

Next, we generalize lines (\ref{lines}).  Given an initial state $(x_0,x_1,x_2)\in\Delta_A$, let us define
$z_0(n):=(x_0,x_1,1-x_0-x_1)\bm P_A\bm P_B^{n}(1,0,0)^\T$ for each $n\ge1$.
Then, for $n$ odd, $z_0(n)<\pi_0$ if and only if
$\alpha_n x_0+\beta_n x_1<\gamma_n$,
and, for $n$ even, $z_0(n)\ge\pi_0$ if and only if
$\alpha_n x_0+\beta_n x_1\le\gamma_n$,
where the coefficients can be expressed in terms of the spectral representation of the matrix $\bm P_B$:
\begin{eqnarray*}
\alpha_n&:=&(-1)^{n-1}(1,0,-1)\bm P_A\bm P_B^{n}(1,0,0)^\T\\
&\;=&(-1)^{n}(3\phi-2)[e_2^{n}(1+\rho^2+S)-e_1^{n}(1+\rho^2-S)]/(4S),\\
\beta_n&:=&(-1)^{n-1}(0,1,-1)\bm P_A\bm P_B^{n}(1,0,0)^\T=(-1)^{n} (3\phi-2)(e_2^{n}-e_1^{n}) (1+\rho^2)/(2S),\\
\gamma_n&:=&(-1)^{n-1}[\bm\pi-(0,0,1)\bm P_A\bm P_B^{n}](1,0,0)^\T\\
&\;=&(-1)^{n}\{e_2^{n} [\phi (1+\rho+\rho^2) (3+3\rho^2+S)-(1+\rho^2) (1+2\rho+3\rho^2+S)]\\
&&\qquad\quad\;{}-e_1^{n} [\phi (1+\rho+\rho^2) (3+3\rho^2-S)-(1+\rho^2) (1+2\rho+3\rho^2-S)]\}\\
&&\qquad/[4 (1+\rho+\rho^2) S].
\end{eqnarray*}
Notice that these definitions are consistent with $(\alpha_1,\beta_1,\gamma_1)$ and $(\alpha_2,\beta_2,\gamma_2)$ defined earlier.

The line $\alpha_n x_0+\beta_n x_1=\gamma_n$ in the plane intersects the line $x_0=\pi_0$ at $x_1=a_n:=(\gamma_n-\alpha_n\pi_0)/\beta_n$;
it intersects the line $x_1=0$ at
$x_0=b_n:=\gamma_n/\alpha_n$;
and it intersects the line $x_0+x_1=1$ at
$x_0=c_n:=(\gamma_n-\beta_n)/(\alpha_n-\beta_n)$.
Each of the lines $\alpha_n x_0+\beta_n x_1=\gamma_n$ $(n\ge1)$ passes through the point $((\phi-2\pi_0)/(3\phi-2),(\phi-2\pi_1)/(3\phi-2))$,
which lies to the left of the line $x_0=\pi_0$ because $\pi_0>1/3$.

We turn next to region 9 of Figure \ref{region-fig}.  This region is determined by $b_2\ge b_1$ (equivalent to $\phi\ge\phi_2$) and $a_1>1-\pi_0$, from which it follows that $b_1<1$ (see Figure \ref{region-fig}).  Consequently, $\Delta_A-\Delta_{ABBA}$ is the triangular region with vertices $(1,0,0)$, $(b_1,0,1-b_1)$, and $(c_1,1-c_1,0)$.  We claim that (a) $\bm P_A\bm P_B$ maps these three points, and hence the triangular region they determine, into $\Delta_{ABBA}$.  Further, we claim that (b) $\bm P_A\bm P_B^2$ maps the four corners of the closure of $\Delta_{ABBA}$, namely $(\pi_0,0,1-\pi_0)$, $(b_1,0,1-b_1)$, $(c_1,1-c_1,0)$, and $(\pi_0,1-\pi_0,0)$, into $\Delta_{ABBA}$, hence maps the closure of $\Delta_{ABBA}$ into $\Delta_{ABBA}$.  This is enough to show that, starting from $\Delta_A$, pattern $ABB$ repeats forever, possibly after an initial $AB$.  To verify (a) and (b), we let
\begin{eqnarray*}
(f_0,f_1,f_2)&:=&(1,0,0)\bm P_A\bm P_B,\\
(g_0,g_1,g_2)&:=&(b_1,0,1-b_1)\bm P_A\bm P_B,\\
(h_0,h_1,h_2)&:=&(c_1,1-c_1,0)\bm P_A\bm P_B,\\
\noalign{\smallskip}
(r_0,r_1,r_2)&:=&(\pi_0,0,1-\pi_0)\bm P_A\bm P_B^2,\\
(s_0,s_1,s_2)&:=&(b_1,0,1-b_1)\bm P_A\bm P_B^2,\\
(t_0,t_1,t_2)&:=&(c_1,1-c_1,0)\bm P_A\bm P_B^2,\\
(u_0,u_1,u_2)&:=&(\pi_0,1-\pi_0,0)\bm P_A\bm P_B^2.
\end{eqnarray*}
First, since $b_2\ge b_1$ in region 9, it follows that if $(x_0,x_1,x_2)\in\Delta$ satisfies $x_0\ge\pi_0$ and $\alpha_1 x_0+\beta_1 x_1<\gamma_1$, then $(x_0,x_1,x_2)\in\Delta_{ABBA}$.  Thus, it suffices to check these two inequalities for each of the seven points listed.  We have done this algebraically but omit the details, as it is reasonably straightforward.

We finally consider regions 10 and 11 of Figure \ref{region-fig}, which can be treated simultaneously.  The regions are determined by $a_1\le1-\pi_0$ (equivalently, $\phi\ge1-\rho/3$), $b_2\ge b_1$ (equivalently, $\phi\ge\phi_2$), and $\alpha_1 f_0+\beta_1 f_1<\gamma_1$, where $(f_0,f_1,f_2):=(1,0,0)\bm P_A\bm P_B$, from which it follows that $a_1\le1-\pi_0$ and $b_1\le1$.  Consequently, $\Delta_A-\Delta_{ABBA}$ is the four-sided region with vertices $(1,0,0)$, $(b_1,0,1-b_1)$,
$(\pi_0,a_1,1-\pi_0-a_1)$, and $(\pi_0,1-\pi_0,0)$, and the closure of $\Delta_{ABBA}$ is the triangular region with vertices $(b_1,0,1-b_1)$, $(\pi_0,a_1,1-\pi_0-a_1)$, and $(\pi_0,0,1-\pi_0)$.  As with region 9, it suffices to show that each of the following seven points $(x_0,x_1,x_2)$ satisfies the inequalities $x_0\ge\pi_0$ and $\alpha_1 x_0+\beta_1 x_1<\gamma_1$ determining $\Delta_{ABBA}$:
\begin{eqnarray*}
(f_0,f_1,f_2)&:=&(1,0,0)\bm P_A\bm P_B,\\
(g_0,g_1,g_2)&:=&(b_1,0,1-b_1)\bm P_A\bm P_B,\\
(h_0,h_1,h_2)&:=&(\pi_0,a_1,1-\pi_0-a_1)\bm P_A\bm P_B,\\
(i_0,i_1,i_2)&:=&(\pi_0,1-\pi_0,0)\bm P_A\bm P_B,\\
\noalign{\smallskip}
(s_0,s_1,s_2)&:=&(b_1,0,1-b_1)\bm P_A\bm P_B^2,\\
(t_0,t_1,t_2)&:=&(\pi_0,a_1,1-\pi_0-a_1)\bm P_A\bm P_B^2,\\
(u_0,u_1,u_2)&:=&(\pi_0,0,1-\pi_0)\bm P_A\bm P_B^2.
\end{eqnarray*}
Again, we have done this algebraically but omit the details. See Figure \ref{regions3,11-fig}.
\end{proof}

Let us describe what appears to happen in each of the eight regions not covered by Theorem \ref{thm-regions3,9--11}.

\begin{conjecture}\label{conj1}
If $(\rho,\phi)$ belongs to region 1 of Figure \ref{region-fig}, then there is a unique limit cycle, which is asymptotically stable and of the form $[1,2]$, as well as an unstable equilibrium.  Indeed,
\begin{eqnarray*}
\Delta_A=\bigcup_{k=0}^\infty\Delta_{(ABBBBABB)^k\overline{ABB}}\cup\bigcup_{k=1}^\infty\Delta_{ABB(ABBBBABB)^k\overline{ABB}},
\end{eqnarray*}
and if the initial state is in $\Delta_B-\Delta_{\overline{B}}$ (see (\ref{Delta(B)})), the trajectory eventually enters $\Delta_A$.

For $(\rho,\phi)$ belonging to region 2 of Figure \ref{region-fig}, there is a unique limit cycle, which is asymptotically stable and of the form $[1,2]$, as well as an unstable equilibrium.  Indeed,
$\Delta_A=\Delta_{\overline{ABB}}\cup\Delta_{ABBBB\overline{ABB}}$,
and for initial states in $\Delta_B-\Delta_{\overline{B}}$, the trajectory eventually enters $\Delta_A$.

If $(\rho,\phi)$ belongs to region 4, 5, 6, 7, or 8 of Figure \ref{region-fig}, then there is a unique limit cycle, which is asymptotically stable and of the form $[1,2]$, as well as an unstable equilibrium.  Indeed,
$$
\Delta_A=\Delta_{A\overline{B}}\cup\bigcup_{k=1}^\infty\Delta_{AB^k\overline{ABB}}.
$$
(Note that $\Delta_{AB^2\overline{ABB}}=\Delta_{\overline{ABB}}$.)
If the initial state is in $\Delta_B-\Delta_{\overline{B}}-\Delta_{BA\overline{B}}-\Delta_{BBA\overline{B}}$, the trajectory eventually enters $\Delta_A-\Delta_{A\overline{B}}$.

If $(\rho,\phi)$ belongs to region 12 of Figure \ref{region-fig}, then there is a unique limit cycle, which is globally asymptotically stable and of the form $[1,2]$.  Indeed,
\begin{eqnarray*}
\Delta_A=\bigcup_{k,l=0}^\infty\Delta_{(AB)^k ABB(AB)^l\overline{ABB}},
\end{eqnarray*}
and only finitely many of the sets comprising the union are nonempty.  If the initial state is in $\Delta_B$, the trajectory eventually enters $\Delta_A$.
\end{conjecture}

The lower boundary of region 1 is the curve $G_{4,2}=0$.  In the unshaded portion of Figure \ref{region-fig} (where $G_{4,2}<0$ and $\phi>2/3$), we have at least one limit cycle, as shown by Theorem \ref{thm-limitcycles}.  See the remark following the statement of the theorem.  Of course, the theorem does not imply that the limit cycles identified there are the only ones, but we conjecture that this is in fact true.

\begin{conjecture}\label{conj2}
Let $n\ge4$ be even.  The curve $b_{n-2}-\pi_0=0$ lies below $H_{n,n-2}=0$ and above $G_{n+2,n}=0$, and the function $b_{n-2}-\pi_0$ is positive above, and negative below, the curve.

If $(\rho,\phi)$ satisfies $G_{n,n-2}<0$ and $E_{n-2}\ge0$, then there are precisely two limit cycles, which are of the forms $[1,n,1,n-2]$ and $[1,n-2]$, as well as an unstable equilibrium.  Indeed,
\begin{eqnarray*}
\Delta_A&=&\Delta_{\overline{AB^nAB^{n-2}}}\cup\bigcup_{k=0}^\infty\Delta_{(AB^nAB^{n-2})^k\overline{AB^{n-2}}}\\
&&\quad{}\cup\Delta_{AB^{n-2}\overline{AB^nAB^{n-2}}}\cup\bigcup_{k=1}^\infty\Delta_{AB^{n-2}(AB^nAB^{n-2})^k\overline{AB^{n-2}}},
\end{eqnarray*}
and if the initial state is in $\Delta_B-\Delta_{\overline{B}}$, the trajectory eventually enters $\Delta_A$.

If $(\rho,\phi)$ satisfies $E_{n-2}<0$ and $E_{n,n-2}\ge0$, then there is a unique limit cycle, which is asymptotically stable and of the form $[1,n,1,n-2]$, as well as an unstable equilibrium.  Indeed,
$\Delta_A=\Delta_{\overline{AB^nAB^{n-2}}}\cup\Delta_{AB^{n-2}\overline{AB^nAB^{n-2}}}\cup\Delta_{AB^n\overline{AB^nAB^{n-2}}}$,
and if the initial state is in $\Delta_B-\Delta_{\overline{B}}$, the trajectory eventually enters $\Delta_A$.

If $(\rho,\phi)$ satisfies $E_{n,n-2}<0$ and $H_{n,n-2}\ge0$, then there are precisely two limit cycles, which are of the forms $[1,n,1,n-2]$ and $[1,n]$, as well as an unstable equilibrium.  Indeed,
\begin{eqnarray*}
\Delta_A&=&\Delta_{\overline{AB^nAB^{n-2}}}\cup\bigcup_{k=0}^\infty\Delta_{(AB^nAB^{n-2})^k\overline{AB^n}}\cup\Delta_{AB^{n-2}\overline{AB^nAB^{n-2}}}
\\&&\quad{}\cup\bigcup_{k=0}^\infty\Delta_{AB^{n-2}(AB^nAB^{n-2})^k\overline{AB^n}}\cup\Delta_{AB^n\overline{AB^nAB^{n-2}}},
\end{eqnarray*}
and if the initial state is in $\Delta_B-\Delta_{\overline{B}}$, the trajectory eventually enters $\Delta_A$.

If $(\rho,\phi)$ satisfies $H_{n,n-2}<0$ and $b_{n-2}-\pi_0\ge0$, then there is a unique limit cycle, which is asymptotically stable and of the form $[1,n]$, as well as an unstable equilibrium.  Indeed,
\begin{eqnarray*}
\Delta_A=\bigcup_{k=0}^\infty\Delta_{(AB^nAB^{n-2})^k\overline{AB^n}}
\cup\bigcup_{k=0}^\infty\Delta_{AB^{n-2}(AB^nAB^{n-2})^k\overline{AB^n}},
\end{eqnarray*}
and if the initial state is in $\Delta_B-\Delta_{\overline{B}}$, the trajectory eventually enters $\Delta_A$.

If $(\rho,\phi)$ satisfies $b_{n-2}-\pi_0<0$ and $G_{n+2,n}\ge0$, then there is a unique limit cycle, which is asymptotically stable and of the form $[1,n]$, as well as an unstable equilibrium.  Indeed,
\begin{eqnarray*}
\Delta_A=\bigcup_{k=0}^\infty\Delta_{(AB^{n+2}AB^n)^k\overline{AB^n}}
\cup\bigcup_{k=1}^\infty\Delta_{AB^n(AB^{n+2}AB^n)^k\overline{AB^n}},
\end{eqnarray*}
and if the initial state is in $\Delta_B-\Delta_{\overline{B}}$, the trajectory eventually enters $\Delta_A$.
\end{conjecture}

\section{Asymptotic cumulative average profit}\label{profit}

In Section \ref{intro} we stated that, if game $B$ is eventually played forever, we have an asymptotically fair game, whereas if the pattern of games is eventually periodic, we have an asymptotically winning game.  Here we try to justify these assertions.  In Section 4 we found that the periodic patterns $[1,n]$ for even $n\ge2$ and $[1,n,1,n-2]$ for even $n\ge4$ can occur.  If our Conjectures \ref{conj1} and \ref{conj2} are correct, then these are the only periodic patterns that can occur.

\begin{theorem}
We denote the asymptotic cumulative average profit per game played by $\mu_{\overline{B}}$ in the situation where game $B$ is eventually played forever, by $\mu_{[1,n]}$ in the case of a limit cycle of the form $[1,n]$ for even $n\ge2$, and by $\mu_{[1,n,1,n-2]}$ in the case of a limit cycle of the form $[1,n,1,n-2]$ for even $n\ge4$.  Then all (Ces\'aro) limits exist, and $\mu_{\overline{B}}=0$, $\mu_{[1,n]}>0$, and $\mu_{[1,n,1,n-2]}>0$.
\end{theorem}

\begin{remark}
This shows that the greedy strategy exhibits the Parrondo effect when $\phi>2/3$ (with some exceptions) but not when $\phi\le2/3$.
\end{remark}

\begin{proof}
In the situation where game $B$ is eventually played forever,
$\mu_{\overline{B}}=\phi\bm\pi\bm\zeta$,
where $\bm\pi$ is the stationary distribution of $\bm P_B$ and $\bm\zeta:=(2p_0-1,2p_1-1,2p_1-1)^\T$ with $p_0=\rho^2/(1+\rho^2)$ and $p_1=1/(1+\rho)$.
We note that, for given $(z_0,z_1,z_2)\in\Delta$,
$(z_0,z_1,z_2)\bm\zeta=z_0(2p_0-1)+(1-z_0)(2p_1-1)=z_0(2p_0-2p_1)+2p_1-1$.
Hence $\bm\pi\bm\zeta=\pi_0(2p_0-2p_1)+2p_1-1=0$,
where the second (algebraic) step will be used again below.  Thus, $\mu_{\overline{B}}=0$.

In the case of a limit cycle of the form $[1,n]$ for even $n\ge2$,
\begin{eqnarray}\label{mu1n}
\mu_{[1,n]}&=&{\phi\over n+1}\bm\pi_{[1,n]}\bm P_A(\bm I+\bm P_B+\cdots+\bm P_B^{n-1})\bm\zeta.
\end{eqnarray}
Now
\begin{eqnarray*}
\bm\pi_{[1,n]}\bm P_A\bm P_B^m\bm\zeta
&=&\bm\pi_{[1,n]}\bm P_A\bm P_B^m\bm (1,0,0)^\T(2p_0-2p_1)+2p_1-1\\
&=&(\pi_0+E_{n,m}/D_n)(2p_0-2p_1)+2p_1-1\\
&=&(E_{n,m}/D_n)(2p_0-2p_1),
\end{eqnarray*}
so (\ref{mu1n}) can be written as
\begin{eqnarray*}
\mu_{[1,n]}={2 \phi\over n+1}(p_0-p_1)\sum_{m=0}^{n-1}{E_{n,m}\over D_n}.
\end{eqnarray*}
By Proposition \ref{propTh7}, $E_{n,m}<0$ for $m=0,1,\ldots,n-1$.
Since $p_0 < p_1$, it follows that $\mu_{[1,n]}>0$.

In the case of a limit cycle of the form $[1,n,1,n-2]$ for $n\ge4$ even,
\begin{eqnarray}\label{mu1n1n-2}
\mu_{[1,n,1,n-2]}&=&{\phi\over2n}\bm\pi_{[1,n,1,n-2]}[\bm P_A(\bm I+\bm P_B+\cdots+\bm P_B^{n-1})\nonumber\\
&&\qquad\qquad\qquad\quad{}+\bm P_A\bm P_B^n\bm P_A(\bm I+\bm P_B+\cdots+\bm P_B^{n-3})]\bm\zeta.\qquad
\end{eqnarray}
Now
\begin{eqnarray*}
\bm\pi_{[1,n,1,n-2]}\bm P_A\bm P_B^m\bm\zeta
&=&\bm\pi_{[1,n,1,n-2]}\bm P_A\bm P_B^m\bm (1,0,0)^\T(2p_0-2p_1)+2p_1-1\\
&=&(\pi_0+G_{n,m}/I_n)(2p_0-2p_1)+2p_1-1 \\
&=&(G_{n,m}/I_n)(2p_0-2p_1)
\end{eqnarray*}
and
\begin{eqnarray*}
&&\bm\pi_{[1,n,1,n-2]}\bm P_A\bm P_B^n\bm P_A\bm P_B^m\bm\zeta\\
&&\quad{}=\bm\pi_{[1,n,1,n-2]}\bm P_A\bm P_B^n\bm P_A\bm P_B^m\bm (1,0,0)^\T(2p_0-2p_1)+2p_1-1\\
&&\quad{}=(\pi_0+H_{n,m}/I_n)(2p_0-2p_1)+2p_1-1\\
&&\quad{}=(H_{n,m}/I_n)(2p_0-2p_1),
\end{eqnarray*}
so (\ref{mu1n1n-2}) can be written as
$$
\mu_{[1,n,1,n-2]}={\phi\over n}(p_0-p_1)\sum_{m=0}^{n-1}{G_{n,m}\over I_n}+{\phi\over n}(p_0-p_1)\sum_{m=0}^{n-3}{H_{n,m}\over I_n}.
$$
By Proposition \ref{propTh7}, $G_{n,m}<0$ for $m=0,1,\ldots,n-1$ and $H_{n,m}<0$ for $m=0,1,\ldots,n-3$.  Since $p_0 < p_1$, it follows that $\mu_{[1,n,1,n-2]}>0$.
\end{proof}

\section*{Acknowledgment}
We thank Derek Abbott for valuable advice. J. Lee was supported by the Basic Science Research Program through the National Research Foundation of Korea (NRF)
funded by the Ministry of Education, Science and Technology (2010-0005364).


\begin{thebibliography}{99}

\bibitem[1]{O}
D.C. Osipovitch, C. Barratt, and P.M. Schwartz, Systems chemistry and Parrondo's paradox: Computational models of thermal cycling, \textit{New J. Chem.} {\bf 33} (2009), pp. 2022--2027.

\bibitem[2]{XPYX}
N.-G. Xie, F.-R. Peng, Y. Ye, and G. Xu, Research on evolution of cooperation among biological system based on Parrondo's paradox game, \textit{J. Anhui Univ. Technol.} {\bf 27} (2010), pp. 167--174.

\bibitem[3]{R}
F.A. Reed, Two-locus epistasis with sexually antagonistic selection: A genetic Parrondo's paradox, \textit{Genetics} {\bf 176} (2007), pp. 1923--1929.

\bibitem[4]{S}
R. Spurgin and M. Tamarkin, Switching investments can be a bad idea when Parrondo's paradox applies, \textit{J. Behav. Finance} {\bf 6} (2005), pp. 15--18.

\bibitem[5]{D}
A. Di Crescenzo, A Parrondo paradox in reliability theory, \textit{Math. Scientist} {\bf 32} (2007), pp. 17--22.

\bibitem[6]{APR}
J. Almeida, D. Peralta-Salas, and M. Romera, Can two chaotic systems give rise to order?, \textit{Physica D} {\bf 200} (2005), pp. 124--132.

\bibitem[7]{APA}
A. Allison, D. Abbott, and C. Pearce, State-space visualisation and fractal properties of Parrondo's games, in {\em Advances in Dynamic Games: Applications to Economics, Finance, Optimization, and Stochastic Control}, A.S. Nowak and K. Szajowski, eds., Annals of the International Society of Dynamic Games 7, Birkh\"auser, Boston, 2002, pp. 613--633.

\bibitem[8]{St}
F. Stjernberg, Parrondo's paradox and epistemology --- when bad things happen to good cognizers (and conversely), in {\em Hommage \`a Wlodek. Philosophical Papers Dedicated
to Wlodek Rabinowicz}, T. R\o nnow-Rasmussen, B. Petersson, J. Josefsson, and D. Egonsson, eds., Department of Philosophy, Lund University, Lund, Sweden. 2007, \url{http://www.fil.lu.se/hommageawlodek/site/abstra.htm}..

\bibitem[9]{FA}
A.P. Flitney, J. Ng, and D. Abbott, Quantum Parrondo's games, \textit{Physica A} {\bf 314} (2002), pp. 35--42.

\bibitem[10]{Py}
R. Pyke, On random walks and diffusions related to Parrondo's games, in {\em Mathematical Statistics and Applications: Festschrift for Constance Van Eeden}, M. Moore, S. Froda, and C. L\'eger, eds., Institute of Mathematical Statistics, Lecture Notes--Monograph Series 42, Beachwood, OH, 2003, pp. 185--216.

\bibitem[11]{DP}
L. Din\'is and J.M.R. Parrondo, Optimal strategies in collective Parrondo games, \textit{Europhys. Lett.} {\bf 63} (2003), pp. 319--325.

\bibitem[12]{HA}
G.P. Harmer and D. Abbott, A review of Parrondo's paradox, \textit{Fluct. Noise Lett.} {\bf 2} (2002), pp. R71--R107.

\bibitem[13]{PD}
J.M.R. Parrondo and L. Din\'is, Brownian motion and gambling: From ratchets to paradoxical games, \textit{Contemp.\ Phys.} {\bf 45} (2004), pp. 147--157.

\bibitem[14]{E}
R.A. Epstein, Parrondo's principle: An overview, in {\em Optimal Play: Mathematical Studies of Games and Gambling}, S.N. Ethier and W.R. Eadington, eds., Institute for the Study of Gambling and Commercial Gaming, University of Nevada, Reno, 2007, pp. 471--492.

\bibitem[15]{A}
D. Abbott, Asymmetry and disorder: A decade of Parrondo's paradox, \textit{Fluct. Noise Lett.} {\bf 9} (2010), pp. 129--156.

\bibitem[16]{VC}
C. Van den Broeck and B. Cleuren, Parrondo games with strategy, in {\em Noise in Complex Systems and Stochastic Dynamics II}, Z. Gingl, J.M. Sancho, L. Schimansky-Geier, and J. Kertesz, eds., Proceedings of the SPIE 5471, SPIE, Bellingham, WA, 2004, pp. 109--118.

\bibitem[17]{EL}
S.N. Ethier and J. Lee, Limit theorems for Parrondo's paradox, \textit{Electron. J. Probab.} {\bf 14} (2009), pp. 1827--1862.

\bibitem[18]{B}
E. Behrends, Stochastic dynamics and Parrondo's paradox, \textit{Physica D} {\bf 237} (2008), pp. 198--206.

\end{thebibliography}
\end{document}